\documentclass[reqno]{amsart}

\pdfoutput=1

\usepackage{graphicx}
\usepackage{amsfonts}
\usepackage{amssymb}
\usepackage{enumitem}
\usepackage{baskervillef}
\usepackage{pdfpages}
\usepackage{tempora}

\usepackage[ukrainian,english]{babel}

\newtheorem*{theorem*}{Theorem}
 
\newtheorem*{lemma*}{Lemma}

\newtheorem*{corollary*}{Corollary}

\newtheorem{innercustomgeneric}{\customgenericname}
\providecommand{\customgenericname}{}
\newcommand{\newcustomtheorem}[2]{%
  \newenvironment{#1}[1]
  {%
   \renewcommand\customgenericname{#2}%
   \renewcommand\theinnercustomgeneric{##1}%
   \innercustomgeneric
  }
  {\endinnercustomgeneric}
}
\newcustomtheorem{customthm}{Theorem}
\newcustomtheorem{customlemma}{Lemma}

\newtheorem*{mainlemma*}{Main Lemma}

\newtheorem*{adianRabin*}{The Adian-Rabin Theorem}

\theoremstyle{definition}
\newcustomtheorem{customremark}{Remark}

\newtheorem*{definition*}{Definition}

\newcommand{\pres}[3]{\textnormal{#1} \langle #2 \mid #3 \rangle}

\begin{document}

\iffalse
%\thispagestyle{empty}
{

	\vfill
	\vspace{3cm}
	{\Large \textsc{Translated into English from the original Ukrainian}\par}
	{\Large \textsc{by} \\ \Large \textsc{C.-F. Nyberg-Brodda}\par}
	\vspace{2.0cm}
	{\Large 2023} \\

}
\clearpage 

\fi

%\thispagestyle{empty}

{
\centering
	
	~
	
	\vspace{-1cm}
	{\scshape\LARGE Six Ukrainian Articles on Semigroup Theory \\ by A. K. Sushkevich \par} 
	\vspace{1cm}
	{\Large \textsc{Translated into English}\par}
	{\Large \textsc{by} \\ \Large \textsc{C.-F. Nyberg-Brodda}\par}
	}
	\vspace{0.5cm}
\begin{center}
\today\\
\rule{\textwidth/2}{1.0pt}
\end{center}	
	
	\vspace{1cm}
%\begin{center}
\noindent{\large\textbf{Translator's Preface}}
%\end{center}

\vspace{0.5cm}

\noindent The present work consists of an English translation of six articles, originally written in Ukrainian, by Anton Kazimirovich Sushkevich (1889--1961). To the best of my knowledge, they constitute all the articles written by Sushkevich in Ukrainian. The bibliographic details for the six articles are as follows; the page numbers stated next to the naming of each article is the page number for the translation in the present work. \\

\

\begin{center}
\makebox[\linewidth]{\hfill \textbf{[Sush1935a]} \hfill \llap{(pp. 6--11)}}
\end{center}
A. K. Sushkevich, \textit{On extending a semigroup to a whole group}, \\ Scientific Notes of the Kharkiv Mathematical Society, \textbf{4}:12 (1935), pp. 81--87. \\

\begin{center}
\makebox[\linewidth]{\hfill \textbf{[Sush1935b]} \hfill \llap{(pp. 12--13)}}
\end{center}
A. K. Sushkevich, \textit{On some properties of a type of generalised groups}, \\ Scientific Notes of the Kharkiv State University, \textbf{2--3} (1935), pp. 23--25 \\

\begin{center}
\makebox[\linewidth]{\hfill \textbf{[Sush1936]} \hfill \llap{(pp. 14--15)}}
\end{center}
A. K. Sushkevich, \textit{Investigations in the field of generalised groups}, \\ Scientific Notes of the Kharkiv State University, \textbf{6--7} (1936), pp. 49--52. \\

\begin{center}
\makebox[\linewidth]{\hfill \textbf{[Sush1937a]} \hfill \llap{(pp. 16--23)}}
\end{center}
A. K. Sushkevich, \textit{On groups of matrices of rank $1$}, \\ J. Inst. Math. Akad. Nauk UkrSSR, \textbf{3} (1937), pp. 83--94. \\

\begin{center}
\makebox[\linewidth]{\hfill \textbf{[Sush1937b]} \hfill \llap{(pp. 24--33)}}
\end{center}
A. K. Sushkevich, \textit{On some types of singular matrices}, \\ Scientific Notes of the Kharkiv State University, \textbf{10} (1951), pp. 5--16. \\

\begin{center}
\makebox[\linewidth]{\hfill \textbf{[Sush1939]} \hfill \llap{(pp. 34--38)}}
\end{center}
A. K. Sushkevich, \textit{Generalised groups of some types of infinite matrices}, \\ Notes Sci. Res. Inst. Khark. State Univ. \& Khark. Math. Soc \textbf{16} (1939), pp. 115--120. \\

\

I first wish to thank C. Hollings (University of Oxford) for providing me with otherwise very inaccessible scans of the Ukrainian originals of the articles, without which this translation could not have been made. 
 
I will attempt to give a summary of each of the articles. This is a non-trivial task; the articles cover rather a broad range of topics in semigroup theory. To read them, a modern semigroup theorist (or indeed anyone familiar with undergraduate-level abstract algebra) will only need to know the following definitions. All italicised words are terms used by Sushkevich; all properties are defined for some fixed semigroup $S$. 
\begin{enumerate}
\item A \textit{generalized group} (or simply \textit{group}!) is a semigroup, i.e. a set equipped with an associative binary operation, hence satisfying what Sushkevich calls the \textit{group property}: the product of two elements in the set again lies in the set.\footnote{Note that this group property is in the spirit of how Galois defined groups -- namely, as a collection of permutations closed under composition.}
\item An \textit{ordinary group} is a group in the modern sense of the word.\footnote{Sushkevich also calls this a \textit{classical} group. Indeed, Hollings notes that while Sushkevich first favoured the term \textit{ordinary},  ``some time later [in 1937], he suggested the name \textit{classical} group'' (see p. 515 and \S5 of Hollings, \textit{Arch. for the Hist. of Exact Sci.} \textbf{63:5}, 2009). However, Hollings does not note that in \textbf{[Sush1939]} Sushkevich had reverted back to \textit{ordinary} in favour of \textit{classical}.}
\item A \textit{semigroup} is a modern-day cancellative semigroup, i.e. a semigroup in which both $ax = ay$ implies $x = y$ (the left cancellative property) and $xa = ya$ implies $x = y$ (the right cancellative property) for all $a, x, y \in S$.
\item The \textit{law of left-sided/right-sided unique reversibility} is simply the right/left cancellative property (note the reversal of the handedness). 
\item The \textit{law of left-sided/right-sided unlimited reversibility} is the law: the equation $Xa = b$ (resp. $aX = b$), for fixed elements $a, b\in S$, has infinitely many solutions in $S$ for the unknown $X$. 
\end{enumerate}

To emphasise point (1), when Sushkevich speaks of \textit{groups}, he means what are today called semigroups (i.e. no element is required to have an inverse). 

As a one-sentence summary, the articles deal with \textit{elementary properties of semigroups}. By \textit{elementary}, I certainly do not mean \textit{easy}. I mean rather those properties which follow from general considerations and from prescribing only, for example, that a certain law (in addition to associativity) should hold. This is particularly clear in \textbf{[Sush1935a]} and \textbf{[Sush1935b]}; in the latter, Sushkevich shows some direct properties of generalized groups in which the right-sided law of unlimited reversibility and the left-sided law of unique reversibility hold. 

The article \textbf{[Sush1935a]} requires particular care. This is because its main result is wrong; Sushkevich claims to have proved, in modern terms, that ``any cancellative semigroup can be embedded in a group''. However, Mal'cev gave a counterexample in 1937 to this statement. We refer, primarily, the reader to Hollings (\textit{Arch. for Hist. of Exact Sci.} \textbf{68}:5 (2014), pp. 641--692) for references and the many unexpected twists in this subject. The error in Sushkevich's argument is quite subtle to extract, and indeed the explanation given by Hollings (which is based on an explanation given by Kurosh) for the error is not quite accurate. In modern terms, Sushkevich's argument goes as follows: let $M$ be a cancellative semigroup, and let $\overline{M}$ be an anti-isomorphic copy of $M$, i.e. if the elements of $M$ are $A, B, \dots$, then the elements of $\overline{M}$ will be $\overline{A}, \overline{B}, \dots$, and if $A \cdot B = C$ in $M$, then $\overline{B} \cdot \overline{A} = \overline{C}$ in $\overline{M}$. Sushkevich then adds an identity element $E$ resp. $\overline{E}$ to each of the semigroups, and forms the monoid free product $M \ast \overline{M}$, i.e. the semigroup free product of the two monoids amalgamating the identity elements. Sushkevich then quotients this free product by the relations $a\overline{a} = \overline{a}a = E$ for all $a \in A$. In this way, he obtains a group, which he calls $G(M)$.\footnote{Here we find Kurosh's misplaced objection: it is not clear to Kurosh that $G(M)$ is a group, or indeed that $G(M)$ satisfies the associative property. However, it is clear that $G(M)$ is a group by using my reformulation via free products -- though not, perhaps, from Sushkevich's original article.} Since all the elements of $M$ can be represented as words over $A$, Sushkevich now claims that $M$ embeds in $G(M)$. But, alas, he never verifies that words representing different elements in $M$ do not become identified in $G(M)$ (and, indeed, he could not have done so, as demonstrated by Mal'cev). 

Indeed, if $M$ is defined by the monoid presentation $\pres{Mon}{A}{u_i = v_i \: (i \in I)}$, i.e. $M$ is generated by some set $A$ subject to the defining relations $u_i = v_i$ for all $i \in I$, then we can let $\overline{A}$ be a set in bijective correspondence with $A$ and such that $A \cap \overline{A} = \varnothing$, and defining
\[
M = \pres{Mon}{\overline{A}}{\overline{u}_i^{\text{rev}} = \overline{v}_i^{\text{rev}} \: (i \in I)},
\]
where $\overline{u}_i$ is the word $u_i$ spelled over the new alphabet $\overline{A}$, and ${}^{\text{rev}}$ denotes word reversal (e.g. $(abacb)^{\text{rev}} = bcaba$). Then $\overline{M}$ is anti-isomorphic to $M$, and Sushkevich's ``group'' will have the presentation: 
\[
G(M) = \pres{Mon}{A}{u_i = v_i, \overline{u}_i^{\text{rev}} = \overline{v}_i^{\text{rev}} \: (i \in I), a\overline{a} = \overline{a}a = 1 \: (\forall a \in A)}.
\]
Of course, if $M$ is group-embeddable, then $M$ certainly embeds in $G(M)$, indeed in this case the natural homomorphism induced by $A \mapsto A$ is injective. But there is no a priori reason to expect that this map is injective in general, and indeed Sushkevich does not demonstrate that it is. 

In spite of its incorrect result, the article \textbf{[Sush1935a]} is remarkable. It shows a clear fluency with semigroup presentations and their relation to group presentations in a way that would not appear in the Soviet literature (or indeed anywhere else) until the work by Markov, Novikov, and Adian over a decade later. He also demonstrates a clear understanding of the relation of his result to other areas of abstract algebra, such as the embedding of integral domains into fields (in \S7). 

The next article, \textbf{[Sush1936]}, contains no new results, but provides a summary of Sushkevich's research carried out while in Kharkiv. This is mostly self-explanatory, but one curious, and to me unexplained, note stands out. On the second page, one reads of a result on real-valued functions proved by ``my student M. Wojdys\l awski''. This is almost certainly the topologist Menachem (Maniek) Wojdys\l awski (1918--1942), a Polish prodigy, who in 1940 studied as a PhD student under Stanis\l aw Mazur. I have found no information about him having ever been a student working with Sushkevich, so this seems like new information, if accurate; I do not, however, know when this studentship would have taken place. Wojdys\l awksi, who was Jewish, was murdered during the Nazi occupation of Poland as part of the destruction of the Cz\c estochowa  Ghetto. For a more or less complete biography, see Maligranda \& Prytu\l a (\textit{Wiad. Mat.} \textbf{49}:1, 2013, pp. 29--66, in Polish). 

The last three articles are all on semigroups of matrices. I will not attempt any lengthy summary, as the results are somewhat technical to state, and the articles themselves straightforward to read. I will, however, mention that they are all infused with the hallmarks of an excellent storyteller, and particularly so \textbf{[Sush1937a]}. This begins by defining a somewhat esoteric-seeming product, and then moves through technical hoops to prove properties about it. Spoiler alert: in the last section, \S7, the truth is revealed. The multiplication defined at the very beginning is nothing but matrix multiplication for a particular class of matrices! Furthermore, certain groups $\mathfrak{G}_{a,b}$ arising in previous sections are shown to be isomorphic to the multiplicative group of all non-zero square roots of the field with which one began -- a very direct definition! Reading the paper again becomes a pleasant experience after this spoiler. Finally, the last article, \textbf{[Sush1939]}, works with rather an unusual subject: infinite matrices. 

In 1997, Wilhelm Magnus (\textit{Amer. Math. Monthly} \textbf{104}:3, p. 269) wrote ``mathematicians understand each other no matter where they come from. [...] Nothing is more international than the community of mathematics''. Sushkevich embodies this spirit to a remarkable extent, having written his articles in five languages: Russian, Ukrainian, English, French, and German. For \textbf{[Sush1936]}, a French translation was already available; I do not doubt that this was written by Sushkevich himself, and it is generally well written (however, it bears mentioning that there are some unfortunate omissions, such as a missing footnote on page 2, omitted sentences, and several references in his list of bibliography missing). There is also a curious typo in \textbf{[Sush1939]} which can probably be directly attributed to his polyglotism; namely, on p. 118 of the Ukrainian original,  at the very beginning of \S3, one reads: ``vcix'', rather than the correct ``\textsc{b}six'' (pronounched \textit{vsikh}, meaning ``all''). This insertion of the Latin letter \textit{v} in the place of its Cyrillic counterpart \textsc{b} is unlikely to be a mere typesetting issue (given how different the letters look), and seems to instead reflect Sushkevich's ability to write in all five languages listed above.

On the subject of typesetting, my general guiding principle has been to maintain, as much as possible, the feel of the original article. My translations are not facsimilies of the originals, but I have preserved all equation numberings, section titles, and notation. I have in some places made very light modifications to the internal structure of the article, e.g. by adding a Theorem environment for clarity. There are some things that make the original article difficult to typeset with perfect accuracy; for example, sometimes $x$ and $\alpha$ appear very similar in the original typesetting, which is exacerbated by the fact that they are both often used as indexing variables in the same sums. It is clear that even the original Ukrainian typesetters have struggled at times, and in many places placed the wrong symbol. I have attempted to fix as many of these issues as I have been able to. Any remaining (or new) errors are, of course, my own. 

\begin{center}
\rule{2cm}{0.5pt}
\end{center}

The original language of these articles is Ukrainian, and were written by a mathematician working in Kharkiv. Ukraine has a long and proud tradition of mathematical excellence, with countless Ukrainian mathematicians making significant contributions to the field. Sushkevich's own mathematical path was beset by the horrors of war and inhumanity, having lived through two world wars and the murder of his former student Wodjys\l awski in the Holocaust (and doubtlessly countless other people known to Sushkevich). Today, almost beyond belief, we find that these horrors are again repeating themselves in Ukraine, following the Russian invasion. This invasion has been accompanied by state-sponsored attempts to erase the culture, heritage, history, and language of the people of Ukraine. It is imperative that we, as mathematicians, continue to promote and preserve the use of Ukrainian as a language of mathematical communication and education, both within Ukraine and internationally, and continue to support our Ukrainian colleagues. With this translation, I hope to contribute to this goal, and make Ukrainian mathematical knowledge more accessible.

Five out of the six original articles in this translation were published in journals directly associated to Kharkiv. Kharkiv, the second largest city of Ukraine, was proclaimed by the Bolsheviks as the capital city of the Ukrainian SSR in 1919; it remained so until 1934. Sushkevich had moved there in 1914, leaving (because of the war) for Voronezh in 1918, and returning only in 1929. Soon thereafter, he was appointed as a full professor at Kharkiv State University. He remained there throughout the Nazi occupation during the Second World War. Kharkiv was the focus of much fighting during the Nazi invasion, and rapidly changed hands multiple time between Nazi Germany and the USSR. Upon its final liberation in 1943, it was ruined; seventy percent of the city had been destroyed, and less than a third of its pre-war population remained. When reading about such horrors, one cannot help but think that it belongs to a bygone age; and indeed, until 2022, the disambiguation page on Wikipedia for ``Battle of Kharkiv'' only listed five such battles, all of them prior to 1943. But today it lists seven. May it never lengthen again. May 2023 bring peace to my Ukrainian friends, colleagues, and the people of Ukraine.

\

\begin{center}
\begin{otherlanguage}{ukrainian}
Слава Україні! Героям слава!
\end{otherlanguage}
\end{center}

\vspace{1cm}

\begin{flushright}
\noindent\textbf{C. F. Nyberg-Brodda}\\
\noindent{Universit\'e Gustave Eiffel}\\
\today
\end{flushright}



\section*{Supplementary material}

# ON EXTENDING A SEMIGROUP TO A WHOLE GROUP

## A. K. Sushkevich  *(Kharkiv)*

**1.** In the sequel, I shall make reference to my work "Über Semigruppen" ("On semigroups"), which was printed in volume VIII of the 4th series of these "Communications". A semigroup differs from a group in that, instead of the law of unlimited reversibility, only the law of unique reversibility is correct for it, that is, not every semigroup has an identity element, and – if there is an identity in the semigroup – not every element has an inverse element in the semigroup. A semigroup proper is necessarily infinite, because in finite systems the law of unique reversibility implies the law of unlimited reversibility. In my aforementioned work, I proved that a semigroup $\mathfrak{S}$ in the most general case can be divided into two parts:

$$\mathfrak{S} = \mathfrak{G} + \mathfrak{H};$$

where $\mathfrak{G}$ (which is called the "group part") and is an ordinary (finite or infinite) group, which may consist of only a single element $E$ (identity element) or may even be empty, while $\mathfrak{H}$ (the "principal part") is a semigroup without identity element, – necessarily infinite, and necessarily present. Furthermore, we have:

$$\text{if } A \subset \mathfrak{G}, \text{ then } \mathfrak{H}\mathfrak{G} \subset \mathfrak{H}, A\mathfrak{H} \subset \mathfrak{H}; \tag{1}$$

$$\text{if } P \subset \mathfrak{H}, \text{ then } \mathfrak{G}P \subset \mathfrak{H}, P\mathfrak{G} \subset \mathfrak{H}; \mathfrak{S}P \subset \mathfrak{H}, P\mathfrak{S} \subset \mathfrak{H}. \tag{2}$$

In the sequel we shall prove that, given a semigroup, it is possible to add a "complement" to the given semigroup to make it into a group, through the addition of new elements (which justifies the term "semigroup").

**2.** Let us first consider the case when the group part is missing, i.e. $\mathfrak{S} = \mathfrak{H}$ is a semigroup without identity; then it is obvious that none of its elements has an inverse. Let $P, Q, R, S, \ldots$ be distinct elements of $\mathfrak{H}$. We create new elements $\overline{P}, \overline{Q}, \overline{R}, \overline{S}, \ldots$ such that each element of $\mathfrak{H}$ corresponds to the the new element denoted by the same letter with a dash above it. We establish the following rule for the composition of these new elements: *whenever $PQ = R$, we set $\overline{Q}\overline{P} = \overline{R}$.* In this way, the set of new elements also forms a semigroup $\overline{\mathfrak{H}}$, which is "anti-isomorphic" with $\mathfrak{H}$.

Next, we create a new element $E$, which we will consider a two-sided identity element for all elements, from $\mathfrak{H}$ as well as from $\overline{\mathfrak{H}}$. Next, we *define*:

$$X\overline{X} = \overline{X}X = E \tag{3}$$

for *every* element $X$ from $\mathfrak{H}$ and the corresponding element $\overline{X} \subset \overline{\mathfrak{H}}$.

Finally, we create several new elements – "products":

$$P\overline{Q}, \overline{P}Q, P\overline{Q}R, \overline{P}Q\overline{R}, P\overline{Q}R\overline{S}, \overline{P}Q\overline{R}S, \ldots \tag{4}$$

such that the "factors" from $\mathfrak{H}$ alternate with the "factors" from $\overline{\mathfrak{H}}$. We will immediately write these "products" without parentheses, so that the associative law is true for them simply by definition.

The "product" of two products (4) is formed by simply writing the factors next to each other, e.g.

$$P\overline{Q}R \cdot U\overline{V} = P\overline{Q}(RU)\overline{V};$$
$$T\overline{S} \cdot \overline{P}Q\overline{R} = T(\overline{SP})Q\overline{R};$$

where we replace $RU$ with a single element using the multiplication table for $\mathfrak{H}$; and, in the same way, $\overline{SP}$ is replaced by a single element using the multiplication table for $\overline{\mathfrak{H}}$.

The products (4) can sometimes be "shortened"; consider, for example

$$P\overline{Q}R\overline{S}T,$$

such that in the multiplication table for $\mathfrak{H}$ we have:

$$R = US; \tag{5}$$

in which case

$$P\overline{Q}R\overline{S}T = P\overline{Q}US\overline{S}T = P\overline{Q}UET = P\overline{Q}(UT) = P\overline{Q}V,$$

if $UT = V$. If, further, in the multiplication table for $\overline{\mathfrak{H}}$ we have

$$\overline{Q} = \overline{W}\overline{V}, \tag{6}$$

then we have:

$$P\overline{Q}R\overline{S}T = P\overline{W}\overline{V}V = P\overline{W}E = P\overline{W}.$$

It should, however, be noted that given $R$ and $S$ there does not always exist a $U$ such that equation (5) holds; and for given $\overline{V}$ and $\overline{Q}$ there does not always exist a $\overline{W}$ such that equation (6) holds, as the law of unlimited reversibility plays no role here whatsoever.

Two products (4) will be considered equal to each other if and only if it is possible to transform one into the other in the specified way, using the equalities (3), (5), and (6). Otherwise, we consider the products (4) to be different.

Let us denote by $\mathfrak{H}_1$ the system of all elements from $\mathfrak{H}$, from $\overline{\mathfrak{H}}$, and all distinct products (4); then $\mathfrak{H}_1$ *is a group*; indeed, the product of any two elements of $\mathfrak{H}_1$ exists and belongs to $\mathfrak{H}_1$; the associative law in $\mathfrak{H}_1$ holds by definition; in $\mathfrak{H}_1$ the element $E$ is an identity element, and for every element from $\mathfrak{H}_1$ there exists an inverse element in $\mathfrak{H}_1$; for example, for

$$P\overline{Q}R\overline{S}T$$

the inverse element is

$$\overline{T}S\overline{R}Q\overline{P},$$

for their product, as it is easy to see, is equal to $E$.

**3.** Let us now consider the general case, when $\mathfrak{S} = \mathfrak{G} + \mathfrak{H}$, i.e. when the group part $\mathfrak{G}$ is also present. We construct for the semigroup $\mathfrak{H}$ the corresponding group $\mathfrak{H}_1$, as in **2**. It is now necessary to combine the groups $\mathfrak{H}_1$ and $\mathfrak{G}$. First of all, we note that the element $E$ defined by (2) is not a new element: we must take $E$ as the identity element of the group $\mathfrak{G}$, as $\mathfrak{G}$ and $\mathfrak{H}_1$ must be subgroups of the group which we want to construct, and subgroups of the same group always have a common identity element. It remains to define the products of elements from $\mathfrak{G}$ with elements from $\overline{\mathfrak{H}}$; we *define*:

$$\begin{cases} A\overline{P} = \overline{Q}, & \text{if } PA^{-1} = Q; \\ \overline{P}A = \overline{R}, & \text{if } A^{-1}P = R; \end{cases} \tag{7}$$

where $A$ is some element from $\mathfrak{G}$ and $P, Q, R$ are elements from $\mathfrak{H}$. Note that then:

$$A = \overline{Q}P = P\overline{R}. \tag{8}$$

In general, if $A$ is any element of $\mathfrak{G}$, and $P$ is any element of $\mathfrak{H}$, then in $\mathfrak{H}$ there are always elements $Q$ and $R$ such that equation (8) holds. Consequently, the entire group $\mathfrak{G}$ is included in the group $\mathfrak{H}_1$ as a subgroup. The relation (8) gives rise to new relations between the products (4): let, for example, $A, B$ be two elements from $\mathfrak{G}$, and for (8):

$$A = \bar{Q}P, \quad B = \bar{S}R;$$

$AB = C$ is also an element of $\mathfrak{G}$, and so by (8):

$$C = \bar{V}U;$$

where $P, Q, R, S, U, V$ are elements from $\mathfrak{H}$; and consequently, we have:

$$\bar{Q}P\bar{S}R = \bar{V}U$$

and there are many more such equalities. That is, with the "inclusion" of the group $\mathfrak{G}$ into the group $\mathfrak{H}_1$, so to speak, it "shrinks": some (even many) elements that were considered different are now combined into one.

**4.** As it known, every group can be defined by a system of "generators"[1] (generating elements) and a system of fundamental relations between the generators: having these, we can build a Cayley table for the group (for an infinite group, this table will, of course, also be infinite; here we can extend and build the table to be arbitrarily large), that is, we will have the entire group with all its laws, including the general laws of our actions. Our point of view is of great importance here: we can immediately see that the basic laws are satisfied, as the associative law, the existence of an identity element, and the existence of inverse elements, are all reflected in the very form of the fundamental relations: when the product of several factors is written without parentheses, it is clear that the associative law is satisfied. Of course, the element $E$ and the inverse elements are all included in the fundamental relations, so they can all be considered to exist. For example, the so-called dihedral group of order 8 (the group of symmetries of a square) is defined as $\{P, Q\}$, or by two generators $P, Q$, between which the following three relations hold:

$$P^4 = E, \quad Q^2 = E, \quad Q^{-1}PQ = P^{-1}.$$

Clearly, in this case we have presupposed the existence of the identity element $E$ and inverse elements. But it is possible, after complicating these fundamental relations and increasing their number, to avoid involving neither $E$ nor inverse elements in them: the relations themselves will already have as their conclusion that both an identity element and inverse elements exist. For the dihedral group of order 8, this results in the following relations:

$$P^5 = P, \quad Q^3 = Q, \quad P^4 = Q^2, \quad QPQ = P^3.$$

Here neither $E$, nor inverse elements, are included explicitly, but e.g. the existence of the identity element $E$ follows from the first three relations, namely $E = P^3 = Q^2$.

In the future, we will consider the associative law to be satisfied in advance. For semi-groups as well as for groups, we will establish a system of generators and fundamental relations holding between them (both the number of generators and relations can be finite or infinite). The difference between a group and a semigroup lies in the fact that for a group the fundamental relations supposes both an identity element as well as the existence of an inverse element for each generator (this is a necessary and sufficient condition for *every element* to have an inverse element); while for a semigroup an identity

---

[1] In the sequel, we will consider them to be independent.

element may exist, but inverses do not exist *for all generators*. If a given semigroup is without a group part, then the fundamental relations do not contain an identity element (and therefore also no inverse elements).

If a given group or semigroup is defined by a system of generators and relations between them, and we take only a part of the entire system of generators, and those relations that only include those generators, then these generators generate a group or a semigroup – a subgroup, or a sub-semigroup, of the given group or semigroup.[2]

**Theorem.** *If a given semigroup has a group part, then the latter is defined, as a group, by those generators of the given semigroup for which the existence of inverse elements follows from the fundamental relations.*

*Proof.* If there is a group part, then there is also an identity element, whose existence consequently follows from the fundamental relations of this semigroup. The group part, like any group, is generated by its generators; we need to prove that apart from the generators which are also generators for the entire semigroup, the group part has no other generators, i.e. – any product of generators of the given semigroup cannot have an inverse element if at least one of the elements of the product does not have one; and this follows directly from Theorem 8 in my previously mentioned work "Über Semigruppen"[3]. □

As for the principal part of a given semigroup, it is generated not only by those generators of this semigroup which do not have inverse elements, because it may happen that the product of some element from the group part by an element from the principal part is such an element, of the principal part, that is not generated only by those generators which have no inverse elements. Thus, in order to form a complete system of generators of the principal part $\mathfrak{H}$ of the given semigroup $\mathfrak{S}$, it is necessary to add a certain (finite or infinite) number of generators to the system of generators for $\mathfrak{S}$, which do not have inverse elements, and which are generated simultaneously by generators from $\mathfrak{H}$ and generators from the group part $\mathfrak{G}$.

5. Let us return to the problem of complementing a given semigroup to turn it into a group. Let $X$ be a given semigroup without a group part (or let it be the principal part of a general semigroup); let $U, V, W, \ldots$ be a complete system of (independent) generators for $\mathfrak{H}$; none of them has an inverse. We *construct* inverses for them: $\overline{U}, \overline{V}, \overline{W}, \ldots$. For each relation

$$UVW \cdots = U_1V_1W_1 \cdots$$

between $U, V, W, \cdots$, we establish the relation

$$\overline{W}\,\overline{V}\,\overline{U} \cdots = \overline{W}_1\overline{V}_1\overline{U}_1 \cdots$$

between $\overline{U}, \overline{V}, \overline{W}$, and construct a semigroup $\overline{\mathfrak{H}}$ with generators $\overline{U}, \overline{V}, \overline{W}, \ldots$, anti-isomorphic with $\mathfrak{H}$; if $UVW \cdots = P$, then we write:

$$\cdots \overline{W}\,\overline{V}\,\overline{U} = \overline{P}.$$

The relation (3) is set to hold only for our generators $U, V, W, \ldots, \overline{U}, \overline{V}, \overline{W}, \ldots$; evidently, it then follows that (3) will hold also for every element $X$ from $\mathfrak{H}$ resp. $\overline{X}$ from $\overline{\mathfrak{H}}$. Next, the construction of the group $\mathfrak{H}_1$ proceeds as in §2 and §3.

---

[2]Moreover, it can happen that a group has a sub-semigroup, or vice versa: a semigroup has a subgroup; for example, the group part of a semigroup is a subgroup of the semigroup in question.

[3]Comm. Khark. Mat. Soc. ser. 4, vol 8 (1934), p. 26.

Let us now take not all, but only some of, the generators $U, V, W, \ldots$ of $\mathfrak{H}$; the chosen generators generate some semigroup $\mathfrak{H}'$ – a sub-semigroup of $\mathfrak{H}$; we complement this semigroup $\mathfrak{H}'$ to a group $\mathfrak{H}'_1$, having first built the semigroup $\overline{\mathfrak{H}'}$, anti-isomorphic with $\mathfrak{H}'$. Next, we consider the products (4) for *all* elements $P, Q, \ldots$ from $\mathfrak{H}$ and for all elements $\overline{R}, \overline{S}, \ldots$ from $\overline{\mathfrak{H}'}$; all these products form a semigroup $\mathfrak{H}_2$, which also contains a subgroup of $\mathfrak{H}'_1$, and the latter is a group part for $\mathfrak{H}_2$. Thus, we can extend the semigroup to this very semigroup, only initially the semigroup did not have a group part, whereas the extended semigroup does. Such an extension is also possible for a semigroup that has a group part. Hence, we will gradually extend this semigroup, finally turning it into a group.

Note that the relations between the generators of the semigroup are *not arbitrary*; the law of unique reversibility leads to great limitations. For example, if these relations imply the existence of a unit for one generator, then they must imply that this unit will also be a unit for all generators of this semigroup.

Note that the generators of a sub-semigroup of a given semigroup are not necessarily only a collection of generators of the given semigroup: they can also be their products. Let $\mathfrak{F}$ be such a sub-semigroup of a given semigroup $\mathfrak{H}$, and let $P = UVW \cdots$ be a generator for $\mathfrak{F}$, where $U, V, W, \ldots$ are generators of $\mathfrak{H}$. By complementing the semigroup $\mathfrak{F}$ to a group in the specified manner, we obtain an inverse element $\overline{P}$ for $P$; but in this way, we also obtain inverse elements for $U, V, W, \ldots$ on the basis of the same Theorem 3 of my work "Über Semigruppen", if the identity element of $\mathfrak{F}$ is also considered the identity element for $\mathfrak{H}$ and the products (4) formed of elements from $\widetilde{\mathfrak{F}}$ and $\mathfrak{H}$. Hence, if all generators from $\mathfrak{H}$ are used in the construction of the semigroup $\mathfrak{F}$, we thereby (that is, with the help of the similarly formed elements $\overline{P}$, the inverses of the generators of the semigroup $\mathfrak{F}$) we complement the entire semigroup $\mathfrak{H}$ into a group; at the same time, of course, it is necessary to establish additional relations required by the law of unique reversibility. For example, if $P = UVW$, and we introduce $\overline{P}$ such that $P\overline{P} = E$, then this implies:

$$UVW\overline{P} = U(VW\overline{P}) = E;$$

hence $VW\overline{P} = \overline{U}$ is the inverse of $U$, which means that the following condition must also be imposed:

$$(VW\overline{P})U = E, \tag{9}$$

because

$$(VW\overline{P})U(VW\overline{P}) = (VW\overline{P})E = E(VW\overline{P});$$

and (9) follows from this by using the law of unique reversibility. Furthermore:

$$V(W\overline{P}U) = E,$$

hence $W\overline{P}U = \overline{V}$ is the inverse of $V$; finally,

$$W\overline{P}UV = W(\overline{P}UW) = E,$$

hence $\overline{P}UV = \overline{W}$ is the inverse of $W$.

*Example:* Let us have a cyclic semigroup:

$$\{A\} = A + A^2 + A^3 + \ldots \text{ in inf.};$$

this has a sub-semigroup:

$$\{A^2\} = A^2 + A^4 + A^6 + \ldots \text{ in inf.};$$

we complement this last one to the whole group; for this, it is necessary to form only a single element $\overline{A^2}$, and a single element $E$, such that:

$$A^2\overline{A^2} = \overline{A^2}A^2 = E;$$

but then $A(A\overline{A^2}) = E$, whence $A\overline{A^2} = \overline{A}$ is the inverse of $A$. Hence, among other relations:

$$\overline{A}A\overline{A^2} = \overline{A}\,\overline{A}$$

$$E\overline{A^2} = (\overline{A})^2, \text{ i.e.: } \overline{A^2} = (\overline{A})^2.$$

**6.** Let us also mention the case when the commutative law holds for the given semigroup; we call such a semigroup *abelian*. Then, when complementing our semigroup to a group, we must assume that every element of the semigroup $\overline{\mathfrak{H}}$ (which is also abelian) commutes with every element of the given semigroup $\mathfrak{H}$. The products (4) have one of two types in this case:

$$R\overline{Q}, \quad \text{or} \quad \overline{Q}R,$$

which are equal to one another; $\mathfrak{H}_1$ is hence abelian.

**7.** The theorem of the extension of a semigroup to a group is quite analogous to the theorem on the extension of an integral domain to a field (the so-called "theorem on forming the field of fractions"); the semigroup here corresponds to the integral domain; the group corresponds to the field; the group part corresponds to the set of so-called "units" of the integral domain; the identity element $E$ corresponds to the so-called "identity element $\varepsilon$". Of course, when proving the theorem on the field of fractions, the method of pairs of elements $(a, b)$ is used, and postulates are made that determine equality, sums, and products, of these pairs.[4] But one can also apply the method we used in §2: an integral domain without an identity element is an abelian semigroup with respect to multiplication; so let us extend it to a group by introducing a single element $\varepsilon$ (if there is none) and for each element $a$ an inverse $\bar{a}$, such that $a\bar{a} = \bar{a}a = \varepsilon$. Next, we introduce new products: $a\bar{b}$. By this very process we have extended our integral domain to a group, with respect to multiplication; it remains to specify the operation of addition for this group; we denote this by the formula:

$$ab + c\bar{d} = (ad + bc)b\bar{d} \qquad (10)$$

with the addition that in the case of a term of the form $\bar{b}$, we write it as: $\varepsilon\bar{b}$. Finally, we set: $0 \cdot \bar{b} = 0$. Now it remains to prove that for this addition, defined by equation (10), all the usual laws of addition are true, that it is connected to the multiplication via the distributive law, and that it coincides with the usual addition for the elements of the given integral domain, which are written as:

$$a = a\bar{\varepsilon}, \quad \bar{\varepsilon} = \varepsilon.$$

The proof of this causes no difficulties.

---

[4] See e.g. van der Waerden, Moderne Algebra, I. Teil, Kap. III, §12. There is a Russian translation of this first part.

# ON SOME PROPERTIES OF A TYPE OF GENERALISED GROUPS

## Prof. A. K. SUSHKEVICH

**1.** This present article is a continuation of my research[1] on an interesting type of generalised infinite groups, defined as those for which the following hold: the associative law, the right-sided law of unlimited reversibility, and the left-sided law of unique reversibility. If $\mathbf{G}$ is such a group and $P$ is any of its elements, then it is easy to see that:

$$P\mathbf{G} = \mathbf{G}; \quad \mathbf{G}P \subset \mathbf{G}.$$

The composition $\mathbf{G}P$ is a group of the same type as $\mathbf{G}$, and consists of all the different elements of the form $GP$ (for all $G \subset \mathbf{G}$)[1], while in the composition $P\mathbf{G}$ every element from $\mathbf{G}$ occurs several times.

**Theorem.** *The group $\mathbf{G}P$ is generally isomorphic with the group $P\mathbf{G}$, i.e. the group $\mathbf{G}$.*

By *generalised* isomorphism I mean, in other words, an *infinite-to-one* isomorphism.

*Proof.* Let $X$ run through all the elements of the group $\mathbf{G}$; let, further, the element $XP$ in $\mathbf{G}P$ correspond to the element $PX$ from $P\mathbf{G}$, i.e. from $\mathbf{G}$; then:

$$\text{if } XP \text{ corresponds to } PX,$$
$$YP \qquad " \qquad PY,$$
$$\text{then } (XPY)P \qquad " \qquad P(XPY),$$

that is, the product $XP.YP$ corresponds to the product $PX.PY$, which also proves our theorem. Let us further note that while for distinct $X, X', X'', \ldots$, the elements $XP, X'P, X''P, \ldots$ may be equal, as the equation $PX = Q$, where $Q$ is a given element of $\mathbf{G}$, has *many* solutions $X, X', X'', \ldots$, for which the products $XP, X'P, X''P, \ldots$ will all be equal. Hence, the obtained isomorphism between $\mathbf{G}P$ and $\mathbf{G}$ is generalised.  □

What is important here is that $\mathbf{G}P$ is a *part* of the group $\mathbf{G}$, and this part is infinite-to-one isomorphic with the entire group $\mathbf{G}$. We can also establish a simple isomorphism between the groups $\mathbf{G}P$ and $\mathbf{G}$, but only by defining a new action in $\mathbf{G}$, as follows:

$$X \circ Y = XPY.$$

Then, when the element $XP$ from $\mathbf{G}P$ corresponds to $X$ from $\mathbf{G}$,
$$\text{and the element } YP \quad " \quad " \qquad " \quad Y \text{ from } \mathbf{G},$$
$$\text{then the product } XP.YP = (XPY)P \text{ corresponds to:}$$
$$XPY = X \circ Y,$$

that is, our new "product" of the elements $X$ and $Y$.

---

[1] See my article "Über einen merkwürdigen Typus der verallgemeinerten unendlichen Gruppen". Comm. of the Kharkiv Math. Soc., ser. 4, vol. 9 (1934), p. 34–44.

**2.** In the theory of ordinary groups, the definition of a group by means of generators (generating elements) and the fundamental relations between them is of great importance. These concepts can also be adapted to generalised groups, both finite and infinite, and in the case of infinite groups, the number of generators can be finite or infinite. All the properties of a given group, as well as the main laws of its action, can be deduced from the relations between the generators of the given group.

**Theorem.** *A group of the type considered here cannot have a finite number of generators.*

*Proof.* Let $A_1, A_2, \ldots A_m$ be a given set of generators: we will prove that they cannot generate a group of our type. Note that every element of the group $(A_1, A_2, \ldots, A_m)$ can be written in the form:

$$P = A_x^{\alpha_x} A_\lambda^{\alpha_\lambda} \cdots A_\rho^{\alpha_\rho}, \tag{1}$$

where the indices $x, \lambda, \ldots, \rho$ (their number is arbitrary, but finite) take values from $1, 2, \ldots, m$ (possibly with repetitions).

It is known[1] that every element of a group of our type has many right identities: let us denote by $\mathbf{E}_1, \mathbf{E}_2, \ldots, \mathbf{E}_m$ the sets of these right identities for our generators $A_1, A_2, \ldots, A_m$ (these sets may have elements in common: it may even happen that one of the sets contains one or more of the others). In addition, we know that the right identities of the elements in the set $E_\lambda$ are all contained in this same set $E_\lambda$. The associative law guarantees that among the identities for the element (1) are all the elements of $\mathbf{E}_\rho$, that is, all the identities of the last generator factor $A_\rho$ of $P$. It follows from this that the set of right identities of any element $P$ of our group contains at least one of these sets $\mathbf{E}_\lambda$.

Consider the element $E_1$ and the set $\mathbf{E}_1$, and denote by $\mathbf{F}_1$ the set of right identities for $\mathbf{E}_1$; as we already know:

$$\mathbf{E}_1 \supset \mathbf{F}_1.$$

The sets $\mathbf{E}_1$ and $\mathbf{F}_1$ cannot coincide with one another, as, for example, $E_1$ belongs to $\mathbf{E}_1$, but does not belong to $\mathbf{F}_1$ (because a group of our type has no idempotent elements). On the other hand, as we have seen, the set $\mathbf{F}_1$ contains at least one of the sets $\mathbf{E}_\lambda$, e.g. $\mathbf{E}_2$, and so:

$$\mathbf{E}_1 \supset \mathbf{F}_1 \supset \mathbf{E}_2.$$

But, reasoning analogously, we can prove that $\mathbf{E}_2$ also contains one of the sets $\mathbf{E}_\lambda$, e.g. $\mathbf{E}_3$, but does not coincide with it; etc. Due to the fact that the number of sets $\mathbf{E}_\lambda$ is finite, after no more than $m$ steps we will reach the conclusion that one of the sets $\mathbf{E}_\lambda$ contains itself as a proper subset; this is a contradiction. In other words, this proves our theorem.        □

---

[1] See my work already quoted above.

# INVESTIGATIONS IN THE FIELD OF GENERALISED GROUPS

## Prof. A. K. SUSHKEVICH

### SUMMARY

During my stay in Kharkiv (since the end of 1929) I have made the following investigations in the field of generalised groups:

I have given a method for representing "left" and "right" groups as well as kernel-type groups[1] by matrices for which the order is higher than the rank. For this, I used rectangular matrices together with Kronecker's theorem that every matrix of order $m$ and rank $n$ (where $m > n$) can be represented as the product of two rectangular matrices: one with $m$ rows and $n$ columns, the other with $n$ rows and $m$ columns, – both of rank $n$. I give a method for transforming a group of ordinary matrices (of order and rank $n$) into a group of matrices of order $m$ and of rank $n$ (with $m > n$), and then build left and right groups from such matrices, and finally, I give a method for constructing such matrices for a given kernel-group; here the problem is reduced to solving a system of equations of the second degree, with several unknowns, with the number of unknowns being greater than the number of variables.

Using an analogous method, I represented, by using rectangular and quadratic matrices, the mixed A. Loewy groups[2] and Brandt groupoids[3].

Next, I tackled the problem of generalising infinite groups in the case when the laws of unlimited and unique reversibility are not fulfilled; but here these laws are substantially different. If we only have the law of unlimited reversibility, dropping the law of unique reversibility, we will obtain what is called a "semigroup", which was known already at the beginning of the 20th century[4]. I studied their structure in detail, and in the general case I proved that they can always be "complemented" into whole groups.

Next, I generalised the notion of a kernel-type group in the case of infinite groups; but here, this type does not have the same meaning as for finite groups: it is only a special case which I have represented by matrices by using a method similar to that which I used for finite generalised groups. This generalisation concerns not only groups, but also semigroups; special cases of this generalisation include "left and right" semigroups and groups. But this still is not the most general case of a "left" (or "right") semigroup, that is an associative generalised group where the left (or the right) law of unique reversibility holds. I have also studied this general case.

Other than this, I have examined a remarkable type of infinite generalised group, of a type which has no analogue in the theory of finite generalised groups; this is an associative group in which the left-hand law of unique reversibility and the law of right-hand

---

unlimited reversibility hold. An example of a group of this type is given by the theory of permutations of a countably infinite set of symbols, by taking those permutations where all the symbols in the lowest row are distinct, but an infinite number of symbols are omitted.[5] My student M. Wojdysławski has found another example of such a group, to wit: the set of all increasing, real-valued functions $f(x)$ with $x$ in the interval $(0, 1)$, and with the conditions $f(0) = 0, f(1) < 1$. The group operation is here, as usual, the substitution of the one function into the other.

Finally, I gave a connection between the "supergroups" ("Übergruppe") of Rauter[6] and ordinary groups, to wit: in addition to the elements of a given group, we also consider all the complexes that can be formed by the elements of the given group, taking into account only the *different* elements which make up the complex; these complexes are multiplied in an ordinary manner. So, all these complexes, together with the individual elements, make up the "supergroup" of Rauter; (I have examined them in detail).

All my research is printed in the following articles of mine.

> *Über die Matrizendarstellung der verallgemeinerten Gruppen.*
> Sci. Notes of the Khark. Math. Soc., ser. 4, vol VI (1933).
>
> *Über Semigruppen.*
> Sci. Notes of the Khark. Math. Soc., ser. 4, vol. VIII (1934).
>
> *Über einen merkwürdigen Typus der verallgemeinerten unendlichen Gruppen*
> Sci. Notes of the Khark. Math. Soc., ser. 4, vol. IX (1934).
>
> *On extending a semigroup to a whole group.*
> Sci. Notes of the Khark. Math. Soc., ser. 4, vol. XII (1935).
>
> *Über eine Verallgemeinerung der Semigruppen.*
> Sci. Notes of the Khark. Math. Soc., ser. 4, vol. XII (1935).
>
> *Über den Zusammenhang des Rauterschen Übergruppe mit den gewöhnlichen Gruppen*
> Math. Zeitschrift, Bd. 38 (1934).
>
> *On some properties of a type of generalised groups.*
> Sci. Notes of Khark. State Univ. For the celebration of the 130th anniversary of Kharkiv State University, 2–3 (1935).

I am currently preparing my monograph on the topic of generalised groups, where all my research mentioned above will be given a more detailed exposition.

The goal of my future research is: find the characteristic features of the type of generalised groups that can be represented by matrices of finite order.

---

[5] This is essentially the same example as given by Baer and F. Levi – "Vollständige irreduzible Systeme von Gruppenaxiomen – Sitzungsber. des Heidelberg Ak., Mathem. nat. Klasse, 1932, 2. Abh.

[6] Rauter. – Abstrakte Kompositionsysteme oder Übergruppen. Journ. Crelle, Bd. 159 (1928).

# ON GROUPS OF MATRICES OF RANK 1

*A. K. Sushkevich* (Kharkiv)

## §1. Definition of elements and their composition (multiplication)

Let $P$ be a given field, – either a number field or abstract; we consider a $n$-dimensional vector space over the field $P$, from which we will take exclusively *unit* vectors, i.e. such vectors $\mathbf{a} = (a_1, a_2, \ldots, a_n)$ such that

$$\mathbf{a}^2 = a_1^2 + a_2^2 + \cdots + a_n^2 = 1. \tag{1}$$

As our *elements* we will take *pairs* of vectors $\mathbf{a}$ and $\mathbf{a}'$ together with a *scalar multiplier* $\alpha$ from $P$. Using uppercase Latin letters, we will write this as follows:

$$A = (\mathbf{a}, \mathbf{a}')\alpha \tag{2}$$

($\mathbf{a}$ is the *left*, and $\mathbf{a}'$ is the *right* vector of the element $A$). We will associate to each element $A$ its *absolute value*, which is determined by the following formula:

$$|A| = \mathbf{a}\mathbf{a}'\alpha. \tag{3}$$

($\mathbf{a}\mathbf{a}'$ is the scalar product of the vectors $\mathbf{a}$ and $\mathbf{a}'$, and is an element in $P$[1]). Among all these elements, we single out those whose scalar multiple is $\alpha = 0$; we consider all such elements to be equal to one another, and simply denote them as 0; this is the *zero element*, or simply *zero*.

Two elements $A = (\mathbf{a}, \mathbf{a}')\alpha$ and $B = (\mathbf{b}, \mathbf{b}')\beta$ which are different from zero (i.e. $\alpha \neq 0, \beta \neq 0$) we will consider as equal if and only if $\mathbf{a} = \mathbf{b}$, $\mathbf{a}' = \mathbf{b}'$, and $\alpha = \beta$.

Now, we define the composition (multiplication) of our elements $A = (\mathbf{a}, \mathbf{a}')\alpha$ and $B = (\mathbf{b}, \mathbf{b}')\beta$ by the following formula:

$$AB = (\mathbf{a}, \mathbf{b}')\alpha\beta\mathbf{a}'\mathbf{b}, \tag{4}$$

where $\mathbf{a}'\mathbf{b}$ is again the usual scalar product of the vectors $\mathbf{a}'$ and $\mathbf{b}$.

From definition (4) it is easy to see that for our multiplication the associative law holds, while the commutative law does not, in general, hold. In particular, from (4) we have:

$$A^2 = AA = (\mathbf{a}, \mathbf{a}')\alpha^2\mathbf{a}\mathbf{a}' = (\mathbf{a}, \mathbf{a}')\alpha|A| = A|A|. \tag{5}$$

Formula (4) implies that $AB = 0$ if one of the elements $A$ or $B$ are equal to 0 (or both $A = 0$ and $B = 0$), or when $\mathbf{a}'b = 0$, i.e. when the vectors $\mathbf{a}'$ and $\mathbf{b}$ are orthogonal; in this last case, we say that the element $A$ is *orthogonal to $B$ on the left*, and $B$ is *orthogonal to $A$ on the right*.

---

[1] If $P$ is a number field, then $|A|$ is not necessarily a positive integer.

## §2. Regular and nilpotent elements

Let us consider the case when for a given element $A$ its absolute value is $|A| = 0$; by **(3)**, this will happen, first when $\alpha = 0$, or further when $A = 0$, or otherwise when $\mathbf{aa}' = 0$, or finally when the left and right vectors of the element $A$ are orthogonal; by (5), in this last case $A^2 = 0$, and we say that $A$ is a *nilpotent element of index* 2. If any element $A$ is orthogonal to $B$ on the left, then as we saw $AB = 0$, and $BA$ is a nilpotent element, because: $BA = (\mathbf{b}, \mathbf{a}')\beta\alpha\mathbf{b}'\mathbf{a}'$, and the vectors $\mathbf{b}$ and $\mathbf{a}'$ are orthogonal.

An element $A$ with absolute value $|A|$ different from zero will be called *regular*. As we have seen, the product of two regular elements can be nilpotent, – precisely when these elements are orthogonal to each other. But the product of two nilpotent elements can in turn be regular – precisely when these elements are not orthogonal to each other. In general, the product of *any* two non-zero elements $A$ and $B$ is regular if and only if $A$ and $B$ are not orthogonal. A nilpotent element is orthogonal to itself.

Let us consider the equation:

$$AX = B \tag{6}$$

where $A$ and $B$ are some given elements. If $A = 0$, and $B \neq 0$, then this equation obviously has no solutions in $X$; if $A = B = 0$, then every element $X$ will be a solution to the equation. So, in the case that $A \neq 0$, and if $A \neq 0$ and $B \neq 0$, then the solutions to the equation (6) will be: the element $X = 0$, and every element $X$ which is orthogonal to $A$ on the right.

Let now $A = (\mathbf{a}, \mathbf{a}')\alpha \neq 0$, $B = (\mathbf{b}, \mathbf{b}')\beta \neq 0$; we denote $X = (\mathbf{x}, \mathbf{x}')\xi$; then from (6) and formula (4) we obtain:

$$(\mathbf{a}, \mathbf{x}')\alpha\xi\mathbf{a}'\mathbf{x} = (\mathbf{b}, \mathbf{b}')\beta. \tag{7}$$

This means that *equation (6) for $A \neq 0$ and $B \neq 0$ has solutions if and only if $\mathbf{a} = \mathbf{b}$, and there are infinitely many such solutions; all of them are given by the formula:*

$$X = (\mathbf{x}, \mathbf{b}')\frac{\beta}{\alpha\mathbf{a}'\mathbf{x}} \tag{8}$$

*where $x$ is an arbitrary (unit) vector.*

Analogously, we have: *the equation*

$$YA = B \tag{6a}$$

for $A \neq 0$ and $B \neq 0$ has solutions if and only if $\mathbf{a}' = \mathbf{b}'$, and there are infinitely many such solutions; all of them are given by the formula:

$$Y = (\mathbf{b}, \mathbf{x}')\frac{\beta}{\alpha\mathbf{a}\mathbf{x}'} \tag{8a}$$

*where $\mathbf{x}'$ is an arbitrary (unit) vector.*

Let us call two elements *left conjugate* if they have the same left vectors; and *right conjugate* if they have the same right vectors; and finally, *generally conjugate* (from both sides) if both the left and the right vectors coincide for them; that is, if the two elements differ only by a scalar multiplier.

The preceding result can now be stated as follows: equation (6) has a solution if and only if $A$ and $B$ are left conjugate; and all solutions $X$ are right conjugate with $B$; for equation (6a) the analogous statement becomes true by interchanging " left" and "right".

Finally, *the equations (6a) and (6a) admit solutions if and only if the elements $A$ and $B$ are generally conjugate; each of them has multiple solutions.*

Note that all of the above results hold regardless of whether the elements $A$ and $B$ are regular, whether they are nilpotent, or whether one of them is regular and the other is nilpotent. Among the solutions $X$ to (6) there will be some that are regular, and some that are nilpotent; in the latter case, the vector $\mathbf{x}$ will be located in the hyperplane perpendicular to the vector $\mathbf{b}'$. The analogous statement also holds for equation (6a).

If the elements $A$ and $B$ are regular and generally conjugate, then the equations (6) and (6a) have one, and only one, solution in common which is also generally conjugate to both $A$ and $B$:

$$X = Y = (\mathbf{a}, \mathbf{a}') \frac{\beta}{\alpha \mathbf{a} \mathbf{a}'} \tag{9}$$

which is also regular. Note that the elements are all generally conjugate; or simultaneously all regular, or simultaneously all nilpotent.

## §3. Idempotent elements

There are elements $E$ such that $E^2 = E$. This includes the element 0. Consider an idempotent element $E = (\mathbf{e}, \mathbf{e}')\varepsilon$ distinct from 0; then

$$E^2 = (\mathbf{e}, \mathbf{e}')\varepsilon^2 \mathbf{e}\mathbf{e}' = E$$

whence $|E| = \mathbf{e}\mathbf{e}'\varepsilon = 1$. Conversely, it is easy to see that an element of absolute value equal to 1 (that is, the multiplicative identity of the field $P$) is idempotent. Thus, we have:

*An element $E$ is idempotent if and only if $|E| = 1$. There are infinitely many idempotent elements; for every regular element $A$ there is a unique idempotent element $I$ generally conjugate to it, which is also the unique two-sided identity element for $A$ (relative to our multiplication for elements) and for all elements generally conjugate to $A$.*

Let $A = (\mathbf{a}, \mathbf{a}')\alpha$; if $I = (\mathbf{a}, \mathbf{a}')\rho$, where $\rho$ is defined by the condition:

$$|I| = \mathbf{a}\mathbf{a}'\rho = 1$$

i.e.

$$\rho = \frac{1}{\mathbf{a}\mathbf{a}'},$$

then it is easy to see that $AI = IA = A$.

It is easy to see that *an idempotent element $E = (\mathbf{e}, \mathbf{e}')\varepsilon$ is a right identity for all elements $X = (\mathbf{x}, \mathbf{e}')\xi$ with which $E$ is right conjugate; and a left identity for all elements $Y = (\mathbf{e}, \mathbf{y})\eta$ with which $E$ is left conjugate, regardless of whether these elements $X, Y$ are regular or nilpotent.*

However, *nilpotent elements do not have a two-sided identity.*

Note also that the scalar multiplier $\varepsilon$ of the idempotent element $E$ is equal to 1 if and only if $\mathbf{e} = \mathbf{e}'$, that is, if $E = (\mathbf{e}, \mathbf{e})$. It is obvious in this case that $E^2 = E$. Conversely, when $\varepsilon = 1$, then we also have $\mathbf{e}\mathbf{e}' = 1$, which means that the unit vectors $\mathbf{e}$ and $\mathbf{e}'$ have one and the same direction, and so they must be equal.

## §4 Generally conjugate elements

Elements of this form, as we have seen, only differ from one another by a scalar factor; they are hence either all regular, or all nilpotent.

*All regular elements generally conjugate to one another form an abelian group, which is isomorphic to the multiplicative group of $P$ (that is, the group with the same multiplication and the same elements, except zero, as $P$).*

Indeed, it is easy to see that the product of two regular, generally conjugate elements is itself a regular element, conjugate with both factors and independent of the order of these factors; the associative law evidently holds; and from §2 it follows that both left- and right-hand unlimited reversibility hold. Let our elements, which are generally conjugate with one another, have the form: $A = (\mathbf{a}, \mathbf{a}')\alpha$; we define the ordinary group that they compose by $\mathfrak{G}_{\mathbf{a}, \mathbf{a}'}$; the identity element here will be the unique idempotent element for them (§3):

$$I = (\mathbf{a}, \mathbf{a}')\rho; \quad \rho = \frac{1}{\mathbf{a}\mathbf{a}'}.$$

The element inverse to $A$ will be:

$$A^{-1} = (\mathbf{a}, \mathbf{a}')\xi, \quad \text{where } \xi = \frac{\rho}{\alpha\mathbf{a}\mathbf{a}'} = \frac{1}{\alpha(\mathbf{a}\mathbf{a}')^2}.$$

Let us now define between the elements $A$ of $\mathfrak{G}_{\mathbf{a}, \mathbf{a}'}$ and the elements $\alpha$ from $P$ (except for zero) the following bijective correspondence:

The element $A(\mathbf{a}, \mathbf{a}')\alpha$ corresponds to $\frac{\alpha}{\rho}$; the element $B = (\mathbf{a}, \mathbf{a}')\beta$ corresponds to $\frac{\beta}{\rho}$; and therefore the element

$$AB = (\mathbf{a}, \mathbf{a}')\alpha\beta\mathbf{a}\mathbf{a}' = (\mathbf{a}\mathbf{a}')\frac{\alpha\beta}{\rho}$$

corresponds to $\frac{\alpha\beta}{\rho^2} = \frac{\alpha}{\rho} \cdot \frac{\beta}{\rho}$; thus, this correspondence is an isomorphism, and a simple one at that, because for different $A$ and $B$, $\frac{\alpha}{\rho}$ and $\frac{\beta}{\rho}$ are also different, and vice versa.

Let us now consider all nilpotent elements which are generally conjugate to one another: $A = (\mathbf{a}, \mathbf{a}')\alpha$ for different $\alpha$ from $P$, all different from zero; such that $\mathbf{a}\mathbf{a}' = 0$. The set of all such elements will be denoted by $\mathfrak{R}_{\mathbf{a}, \mathbf{a}'}$. Any two elements of $\mathfrak{R}_{\mathbf{a}, \mathbf{a}'}$ will always be orthogonal, i.e. their product is always equal to zero, or:

$$\mathfrak{R}_{\mathbf{a}, \mathbf{a}'}^2 = 0.$$

Adding to $\mathfrak{R}_{\mathbf{a}, \mathbf{a}'}$ this element 0, we obtain a so-called " zero group", i.e. such a generalized associative group in which the product of any pair of elements is always equal to 0 – the "kernel" of the group.

## §5 Elements which are conjugate on one side

We will consider the set of all elements with the same vector, i.e. which are right conjugate to one another. This set $\mathfrak{A}_a$ is divided into the complexes (groups) $\mathfrak{G}_{x,a}$ (where $xa \neq 0$) and the complexes $\mathfrak{R}_{y,a}$ (where $ya = 0$); no two of these complexes have elements in common. We can write this symbolically as:

$$\mathfrak{A}_a = \sum_x \mathfrak{G}_{x,a} + \sum_y \mathfrak{R}_{y,a}. \tag{10}$$

(Here, as in the theory of groups, by + and the "sum", we mean not an action, but simply the combination of the individual elements or complexes into a single complex). Every element from $\mathfrak{R}_{y,a}$ is orthogonal on the right with each element from $\mathfrak{G}_{x,a}$; that is (§§1, 2):

$$\mathfrak{G}_{x,a}\mathfrak{R}_{y,a} = 0$$

while

$$\mathfrak{R}_{y,a}\mathfrak{G}_{x,a} = \mathfrak{R}_{y,a}. \tag{11}$$

We also have

$$\mathfrak{G}_{x,a}\mathfrak{G}_{x',a} = \mathfrak{G}_{x,a} \tag{12}$$

because the product of two regular, right conjugate elements is also regular, which is easy to see.

*All right conjugate, regular elements constitute a left group:*

$$\mathfrak{A}'_a = \sum_x \mathfrak{G}_{x,a}. \tag{13}$$

(A left group is an associative group, where the laws of left-hand unlimited and left-hand unique reversibility hold).

Indeed, by (12), the group property for $\mathfrak{A}'_a$ is fulfilled; it is still necessary to prove that equation (6a) for all elements $A$ and $B$ from $\mathfrak{A}'_a$ will always admit exactly one solution $Y \subset \mathfrak{A}'_a$, but this follows directly from §2: if $A = (r,a)\alpha$, $B = (s,a)\beta$, then by equation (8a) from §2, every solution to equation (6a) is of the form:

$$Y = (s,x)\frac{\beta}{\alpha r x}$$

and among these solutions, only *one* belongs to $\mathfrak{A}'_a$, namely for $x = a$: $Y_a = (s,a)\frac{\beta}{\alpha r a}$. It follows that $Y_a$ is generally conjugate with $B$, i.e. it belongs to the same ordinary group $\mathfrak{G}_{s,a}$ as $B$.

Note that equation (6a) always has the solution $Y$ when $A \subset \mathfrak{A}'_a$ or when $B = (y,a)\beta \subset \mathfrak{R}_{y,a}$; but only one of these solutions belongs to $\mathfrak{A}_a$, namely:

$$Y = (y,a)\frac{\beta}{\alpha r a}$$

and this solution belongs to the complex $\mathfrak{R}_{y,a}$ – i.e. the same as $B$ belongs to.

The second part of the set $\mathfrak{A}_a$, i.e. $\sum_u \mathfrak{R}_{y,a} = \mathfrak{R}_a$, together with the element 0, forms a "zero" group, because the product of any pair of its elements is equal to zero.

Therefore, the totality of the set $\mathfrak{A}_a$ together with the element 0 forms, as one might say, a *generalized left group*; its structure is given by the following schema:

|                    | 0 | $\mathfrak{R}_{y,a}$ | $\mathfrak{G}_{x,a}$ |
|--------------------|---|----------------------|----------------------|
| 0                  | 0 | 0                    | 0                    |
| $\mathfrak{R}_{y,a}$ | 0 | 0                    | $\mathfrak{R}_{y,a}$ |
| $\mathfrak{G}_{x,a}$ | 0 | 0                    | $\mathfrak{G}_{x,a}$ |

Table 1

From Table 1, it follows, among other things, that taking not all, but only *some* (arbitrarily many, or an infinite number if $P$ is infinite) of the complexes $\mathfrak{R}_{y,a}$ together with some of the groups $\mathfrak{G}_{x,a}$ and adding the element 0, we also obtain a generalized left group of the same type as the whole group $\mathfrak{A}_a + 0$ with the same schema as in Table 1.

Everything that has, at this point, been said about the generalized left groups can be repeated for the *generalized right group* consisting of all left conjugate elements, i.e. to the same left vector $a$ – together with the element 0. The set $\mathfrak{B}_a$ of such elements can be given in the form:

$$\mathfrak{B}_a = \sum_x \mathfrak{G}_{a,x} + \sum_y \mathfrak{R}_{a,y}. \tag{10a}$$

All $\mathfrak{G}_{a,x}$ are ordinary groups of elements generally conjugate with one another; and $\mathfrak{R}_{a,y} + 0$ is a zero group; no two of the complexes of groups $\mathfrak{G}_{a,x}$ or $\mathfrak{R}_{a,y}$ have any elements in common. We have:

$$\mathfrak{R}_{a,y}\mathfrak{G}_{a,x} = 0; \quad \mathfrak{R}_{a,y'}\mathfrak{R}_{a,y} = 0; \quad \mathfrak{G}_{a,x}\mathfrak{R}_{a,y} = \mathfrak{R}_{a,y} \tag{11a}$$

$$\mathfrak{G}_{a,x}\mathfrak{G}_{a,x'} = \mathfrak{G}_{a,x'}. \tag{12a}$$

*All left conjugate, regular elements constitute a right group:*

$$\mathfrak{B}'_a = \sum_x \mathfrak{G}_{a,x}. \tag{13a}$$

We may also take not all, but only some of the groups $\mathfrak{G}_{a,x}$ and the complexes $\mathfrak{R}_{a,y}$, together with the element 0; we obtain a generalized right group with the same structure as $\mathfrak{B}_a + 0$; and this structure is given in the following schema:

|               | 0 | $\mathfrak{R}_{a,y}$ | $\mathfrak{G}_{a,x}$ |
|---------------|---|----------------------|----------------------|
| 0             | 0 | 0                    | 0                    |
| $\mathfrak{R}_{a,y}$ | 0 | 0             | 0                    |
| $\mathfrak{G}_{a,x}$ | 0 | $\mathfrak{R}_{a,y}$ | $\mathfrak{G}_{a,x}$ |

TABLE 2

## §6 The generalized group of all elements

Denoting the totality of all our elements (for a given field $P$ and for a given $n$) by $\mathfrak{H}$, we can say that $\mathfrak{H}$ is a generalized associative group without the laws of unlimited and unique reversibility[2] with a zero element 0, which is also a kernel. In addition to the element 0, this group $\mathfrak{H}$ consists of all possible complexes $\mathfrak{A}_a$ for all different $a$, or all complexes $\mathfrak{B}_b$ for all different $b$, or, finally, all possible ordinary groups $\mathfrak{G}_{a,b}$ and the complexes $\mathfrak{R}_{a',b'}$ for all different $a, b$ with $ab \neq 0$ and all different $a', b'$ with $a'b' = 0$ (see §5); symbolically, we write this as:

$$\mathfrak{H} = 0 + \sum_a \mathfrak{A}_a = 0 + \sum_b \mathfrak{B}_b = 0 + \sum_{a,b} \mathfrak{G}_{a,b} + \sum_{a',b'} \mathfrak{R}_{a',b'}; \quad (ab \neq 0, a'b' = 0). \tag{14}$$

From formula (4) in §1 we deduce from $a'b \neq 0$ that:

$$\mathfrak{G}_{a,a'}\mathfrak{G}_{b,b'} = \mathfrak{G}_{a,b'} \quad \text{or} \quad = \mathfrak{R}_{a,b'}, \tag{15}$$

depending on whether $ab' \neq 0$ or $ab' = 0$. Under the same condition: either $\mathfrak{R}_{a,a'}\mathfrak{R}_{b,b'}$ (when $aa' = bb' = 0$) or $\mathfrak{G}_{a,a'}\mathfrak{R}_{b,b'}$ (when $aa' \neq 0 = bb'$), or $\mathfrak{R}_{a,a'}\mathfrak{G}_{b,b'}$ (when $aa' = 0 \neq bb'$), or $\mathfrak{R}_{a,a'}\mathfrak{G}_{b,b'}$ (when $aa' = 0 \neq bb'$) are equal to either $\mathfrak{G}_{a,a'}$ (when $ab' \neq 0$), or $\mathfrak{R}_{a,b'}$ (when $ab' = 0$). From this it easily follows that:

$$\mathfrak{A}_a\mathfrak{A}_b = 0 + \mathfrak{A}_b; \quad \mathfrak{B}_a\mathfrak{B}_b = 0 + \mathfrak{B}_a. \tag{16}$$

Furthermore:

$$\mathfrak{B}_a\mathfrak{A}_b = \mathfrak{G}_{a,b} + 0, \quad \text{or} \quad = \mathfrak{R}_{a,b} + 0, \tag{17}$$

_____

[2]For generalized groups, see my works: "Über endliche Gruppen ohne das Gesetz der eindeutigen Umkehrbarkeit", Math. Ann., **99** (1928); "Über eine Verallgemeinerung der Semigruppen", Notes of the Khark. Math. Soc., ser. 4, **12** (1935); see also my monograph "Theory of generalized groups", DNTVU, 1937, chapters III and IV.

depending on whether $ab \neq 0$ or $ab = 0$;

$$\mathfrak{A}_b \mathfrak{B}_a = \mathfrak{H} - 0, \quad \text{or } = 0, \tag{18}$$

depending on whether $ab \neq 0$ or $ab = 0$.

In this way, the group $\mathfrak{H}$ is a collection of groups of "kernel" type, which I have investigated in the works mentioned above. In this group $\mathfrak{H}$, it is possible to distinguish, as subgroups, many groups of "kernel" type. Let $a_1, a_2, \ldots, a_r$ and $b_1, b_2, \ldots, b_s$ be two systems of unit vectors, satisfying that $a_x b_\lambda \neq 0$ ($k = 1, 2, \ldots, r; \lambda = 1, 2, \ldots, s$). We denote by $\mathfrak{G}_{a_x, b_\lambda}$ the following group of elements: $(a_x, b_\lambda)_\alpha$, where $\alpha$ traverses all non-zero elements of $P$. We have $rs$ such groups $\mathfrak{G}_{a_x, b_\lambda}$, which together constitute a group $\mathfrak{K}$ of "kernel" type; here $\mathfrak{A}_{b_\lambda} = \sum_{x=1}^{r} \mathfrak{G}_{a_x, b_\lambda}$ is a left group, and $\mathfrak{B}_{a_x} = \sum_{\lambda=1}^{s} \mathfrak{G}_{a_x, b_\lambda}$ is a right group, and

$$\mathfrak{K} = \sum_{\lambda=1}^{s} \mathfrak{A}_{b_\lambda} = \sum_{x=1}^{r} \mathfrak{B}_{a_x} = \sum_{x=1}^{r} \sum_{\lambda=1}^{s} \mathfrak{G}_{a_x, b_\lambda}.$$

The numbers $r$ and $s$ can also be equal to $\infty$ (if the field $P$ is infinite).

## §7 Matrices of rank 1

Our "elements", which were considered in the previous paragraphs, are nothing but matrices of the $n$th order and rank 1 (except the element 0, which is the zero matrix). Indeed, in a matrix of rank 1, all rows are pairwise linearly dependent, and all columns are also pairwise linearly dependent; hence a matrix of rank 1 can be written in the form:

$$\begin{pmatrix} c_1 c_1' & c_1 c_2' & \cdots & c_1 c_n' \\ c_2 c_1' & c_2 c_2' & \cdots & c_2 c_n' \\ \cdots & \cdots & \cdots & \cdots \\ c_n c_1' & c_n c_2' & \cdots & c_n c_n' \end{pmatrix} \tag{19}$$

where not all $c_x$ and not all $c_\lambda$ are equal to zero, that is:

$$c_1^2 + c_2^2 + \cdots + c_n^2 = k \neq 0, c_1'^2 + c_2'^2 + \cdots + c_n'^2 = l \neq 0.$$

If we set $a_x = \frac{c_x}{\sqrt{k}}$, and $a_\lambda' = \frac{c_\lambda'}{\sqrt{l}}$, we can rewrite (19) as:

$$\sqrt{kl} \begin{pmatrix} a_1 a_1' & a_1 a_2' & \cdots & a_1 a_n' \\ a_2 a_1' & a_2 a_2' & \cdots & a_2 a_n' \\ \cdots & \cdots & \cdots & \cdots \\ a_n a_1' & a_n a_2' & \cdots & a_n a_n' \end{pmatrix} = \sqrt{kl} \begin{pmatrix} a_1 \\ a_2 \\ \cdots \\ a_n \end{pmatrix} (a_1', a_2' \cdots, a_n') \tag{20}$$

here $a_1^2 + a_2^2 + \cdots + a_n^2 = 1, a_1'^2 + a_2'^2 + \cdots + a_n'^2 = 1$, so the rectangular matrices $\begin{pmatrix} a_1 \\ a_2 \\ \cdots \\ a_n \end{pmatrix}$

and $(a_1', a_2' \cdots, a_n')$ can be considered as unit vectors $\mathbf{a}, \mathbf{a}'$, and equation (20) gives us a matrix (19) of rank 1 in the same way that we denoted our "element" by equation (2)

in §1, for $\alpha = \sqrt{kl}$. By the usual rule for matrix multiplication, we have:

$$\begin{pmatrix} a_1a_1' & a_1a_2' & \cdots & a_1a_n' \\ a_2a_1' & a_2a_2' & \cdots & a_2a_n' \\ \cdots & \cdots & \cdots & \cdots \\ a_na_1' & a_na_2' & \cdots & a_na_n' \end{pmatrix} \begin{pmatrix} b_1b_1' & b_1b_2' & \cdots & b_1b_n' \\ b_2b_1' & b_2b_2' & \cdots & b_2b_n' \\ \cdots & \cdots & \cdots & \cdots \\ b_nb_1' & b_nb_2' & \cdots & b_nb_n' \end{pmatrix} =$$

$$= \begin{pmatrix} a_1b_1' & a_1b_2' & \cdots & a_1b_n' \\ a_2b_1' & a_2b_2' & \cdots & a_2b_n' \\ \cdots & \cdots & \cdots & \cdots \\ a_nb_1' & a_nb_2' & \cdots & a_nb_n' \end{pmatrix} \cdot (a_1'b_1 + a_2'b_2 + \cdots + a_n'b_n),$$

which says that the multiplication for our "elements", as defined by equation (4) in §1, is nothing but the usual multiplication of matrices. Therefore, the group $\mathfrak{H}$ considered in §6, and the subgroups considered in the previous paragraphs, are nothing but groups of matrices of order $n$ and rank 1 (together with the zero matrix with rank 0), and $\mathfrak{H}$ is the most general such group.

One more point remains to be discussed. In §§1–6 we took for the scalar factors $\alpha, \beta, \xi, \ldots$ the elements of a field $P$. But it follows from equation (20) that the factor $\alpha = \sqrt{kl}$ can also be irrational, i.e. that it is necessary to extend the field $P$. If $P$ is the field of all real numbers, or more generally some abstract field which contains the square root $\sqrt{\alpha}$ of all its elements $\alpha$, then it is not necessary to extend $P$. But even in the general case, it is not necessary to extend $P$ as a field, via a field extension: we extend only $P$ (without zero) as a *multiplicative group* by considering all ("positive") square roots of all non-zero elements of $P$, or, in other words, we will consider the multiplicative group $P'$ of square roots of non-zero elements of $P$ (it is easy to verify that $P'$ is indeed an ordinary abelian group). Now, let us consider the elements of the vector (2), where $\mathbf{a}, \mathbf{a}'$ are unit vectors over the field $P$, and the scalar factor $\alpha$ is an element of the group $P'$. For such elements, all further conclusions in §§1–6 remain valid; the group $\mathfrak{G}_{a,b}$ is now simply isomorphic to the group $P'$.

# ON SOME TYPES OF SINGULAR MATRICES

*Prof. A. K. Sushkevich*

**§1.** In this article, we will consider the so-called "singular" square matrices of order $m$, that is, matrices whose rank $r$ is less than their order: $r < m$. They are characterised by the fact that they have 0 for one of their characteristic roots, and if this root 0 is a multiplicity of $\alpha$, then $\alpha \geq m - r$, and such a matrix has exactly $m - r$ elementary divisors corresponding to this root 0. In particular, for such a matrix $A$ all elementary divisors corresponding to the root 0 are linear; we will call such a matrix $A$ *regular*; this has the property that all its powers $A^2, A^3, \ldots$ have the same rank $r$; and conversely, this property characterises regular matrices. Finally, one can give yet another definition for regular matrices, as I did in my work "Über die Matrizendarstellung der verallgemeinerten Gruppen" (Notes of the Kharkiv Math. Soc., ser 4, vol. 6, 1933, pp. 27–38): I mention there Kronecker's theorem, that every matrix $A$ of order $m$ and rank $r$ can be presented as a product

$$A = KL \tag{1}$$

of an $mr$-matrix[1] $K$ and an $rm$-matrix $L$, both of rank $r$. If we now consider the product $A_1 = LK$, we see that it is a square matrix of order $r$ and rank $\leq r$. So the matrix $A$ is regular if and only if the matrix $A_1$ is non-singular, i.e. when the rank $= r$. This is the third definition of a regular matrix.

In particular, every idempotent matrix $I$ is regular, precisely because $I^2 = I$. If $I = KL$ is a factorisation analogous to (1) (where we assume $I$ has order $m$ and rank $r$), then it is easy to see that the matrix: $LK = E_r$ is a unit matrix of order $r$. Indeed, the equality $I^2 = I$ gives:

$$KLKL = KE_rL = KL;$$

and so: $LKE_rLK = LKLK$, whence: $E_r^3 = E_r^2$, and due to the fact that $E_r$ is non-singular, it follows from this that

$$E_r^2 = E_r;$$

and therefore $E_r$ is a unit matrix of order $r$, because a unit matrix is the unique non-singular idempotent matrix of a given order.

**Theorem.** *There are many matrices of order $m$ and rank $r$ that have a given idempotent matrix $I$ for its two-sided identity; all these matrices are regular and form an ordinary group.*

*Proof.* Let $A_1$ be a non-singular matrix of order $r$; then $A = KA_1L$ is a matrix of order $m$ and rank $r$,[2] for which there is a two-sided identity; indeed:

$$AI = KA_1LKL = KA_1E_rL = KA_1L = A;$$

analogously, we also have $IA = A$. If $B_1 \neq A_1$ is another non-singular matrix of order $r$, then $B = KB_1L \neq A$, because from $B = A$ it follows that $LBK = LAK$, so

---

[1] That is, a rectangular matrix with $m$ rows and $r$ columns.

[2] See my previously referenced article "Über die Matrizendarstellung der verallgemeinerten Gruppen", §1, Theorem 1.

$B_1 = A_1$; but $B$ also has $I$ as its two-sided identity. Thus the first part of the theorem is proved. Conversely, suppose that $A$ is any matrix of order $n$ and rank $r$, which has a two-sided identity $I$. Let $LAK = A_1$; the matrix $A_1$ has order $r$; we need to prove that its rank $= r$, but this follows from the equality:

$$KA_1L = KLAKL = IAI = A,$$

because the rank of $A$ is $r$. In order to complete the proof of the theorem, it is necessary to prove that such a matrix $A$ has an "inverse" matrix with respect to $I$, i.e. such a matrix $A'$, also of order $m$ and rank $r$, with a two-sided identity $I$, such that

$$AA' = A'A = I.$$

But such a matrix is: $A' = KA_1^{-1}L$, where $A_1^{-1}$ is the inverse matrix of the non-singular matrix $A_1$; $A_1^{-1}$ is also a non-singular matrix of of order $r$. Indeed:

$$AA' = KA_1LKA_1^{-1}L = KA_1E_rA_1^{-1}L = KA_1A_1^{-1}L = KE_rL = KL = I;$$

and analogously for $A'A = I$.[3]                                          □

**Remark.** It is easy to also prove that there are many matrices of order $m$ and rank $p$ that have a given idempotent only as a right (or only as a left) identity. But these matrices do not form a group (even a generalized one), because the group property is not fulfilled for them: among them, there are also non-regular matrices.

**§2.** Every idempotent matrix of order $m$ and rank $r$ has characteristic roots equal to only 0 or 1, and only linear elementary divisors. And from this it follows that every such idempotent matrix has a characteristic root of 1 of multiplicity $r$, and a root 0 of multiplicity $m - r$, and is similar to the following matrix:

$$E_m^{(r)} = \begin{bmatrix} 1 & 0 & \cdots & 0 & 0 & \cdots & 0 \\ 0 & 1 & \cdots & 0 & 0 & \cdots & 0 \\ \cdots & \cdots & \cdots & \cdots & \cdots & \cdots & \cdots \\ 0 & 0 & \cdots & 1 & 0 & \cdots & 0 \\ 0 & 0 & \cdots & 0 & 0 & \cdots & 0 \\ \cdots & \cdots & \cdots & \cdots & \cdots & \cdots & \cdots \\ 0 & 0 & \cdots & 0 & 0 & \cdots & 0 \end{bmatrix} \begin{array}{l} \left.\vphantom{\begin{matrix}1\\0\\ \cdots \\0\end{matrix}}\right\} r \\ \left.\vphantom{\begin{matrix}0\\ \cdots \\0\end{matrix}}\right\} m-r \end{array}$$

$$\underbrace{\phantom{xxxxxxxx}}_{r} \quad \underbrace{\phantom{xxxxxxxx}}_{m-r}$$

because this matrix $E_m^{(r)}$ is nothing but the Jordan normal form for the matrix $I$. Therefore, there exists a non-singular matrix $P$ such that $P^{-1}IP = E_m^{(r)}$.

**Theorem.** *The matrices $P$ and $P^{-1}$ can be obtained by properly "extending" the matrices $K$ and $L$ (where $I = KL$, as in §1) into non-singular, mutually inverse matrices of order $m$.*

*Proof.* The theorem reduces to the following: let us suppose we are given a normalized biorthogonal system consisting of two rows of $r$ $m$-dimensional linearly independent vectors: $\mathbf{a}_\lambda = (a_{\lambda 1}, a_{\lambda 2}, \ldots, a_{\lambda m})$ and $\mathbf{b}_\lambda = (b_{\lambda 1}, b_{\lambda 2}, \ldots, b_{\lambda m})$, for $\lambda = 1, 2, \ldots, r$. Then it is necessary to prove that this biorthogonal system can be "extended" by adding vectors $\mathbf{x}_\mu = (x_{\mu 1}, x_{\mu 2}, \ldots, x_{\mu m})$ to the series of vectors $\mathbf{a}_\lambda$, and to add vectors $\mathbf{y}_\mu = (y_{\mu 1}, y_{\mu 2}, \ldots, y_{\mu m})$ (in both cases, $\mu = 1, 2, \ldots, s$, where $s = m - r$), such that the extended system remains a normalized biorthogonal system.

---

[3]This last part of our theorem provides a refinement of Theorem 4 in §1 of my previously referenced article "Über die Matrizendarstellung...".

So, we are given the condition:

$$\sum_{\lambda=1}^{m} a_{\alpha\lambda} b_{\beta\lambda} = \delta_{\alpha\beta} = \begin{cases} 1 & \text{if } \alpha = \beta, \\ 0 & \text{if } \alpha \neq \beta. \end{cases} \tag{2}$$

where $\alpha, \beta = 1, 2, \ldots, r$. To determine the $x_{\mu\nu}$, we have the following equalities:

$$\sum_{\mu=1}^{m} b_{\alpha\mu} x_{x\mu} = 0 \qquad (\alpha = 1, 2, \ldots, r; \quad x = 1, 2, \ldots, s). \tag{3}$$

For a given $x$, this is a system of $r$ homogeneous linear equations with $m$ unknowns $x_{x1}, x_{x2}, \ldots, x_{xm}$. The rank of the matrix of coefficients of this system is $r$, because the $r$ vectors $\mathbf{b}_\lambda$ are linearly independent. Such a system (3) has $m - r = s$ linearly independent solutions, which we will take as the components of the vectors $\mathbf{x}_1, \mathbf{x}_2, \ldots, \mathbf{x}_s$. Let us now prove that the $r + s = m$ vectors $\mathbf{a}_1, \ldots, \mathbf{a}_r, \mathbf{x}_1, \ldots, \mathbf{x}_s$ are all linearly independent. Suppose:

$$c_1 \mathbf{a}_1 + \cdots + c_r \mathbf{a}_r + d_1 \mathbf{x}_1 + \cdots + d_s \mathbf{x}_s = 0;$$

which yields:

$$\sum_{\alpha=1}^{r} c_\alpha a_{\alpha\lambda} + \sum_{x=1}^{s} d_x x_{x\lambda} = 0 \quad (\lambda = 1, 2, \ldots, m). \tag{4}$$

We multiply this by $b_{\mu\lambda}$ and sum it for $\lambda$ from 1 to $m$:

$$\sum_{\lambda=1}^{m} \sum_{\alpha=1}^{r} c_\alpha a_{\alpha\lambda} b_{\mu\lambda} + \sum_{\lambda=1}^{m} \sum_{x=1}^{s} d_x x_{x\lambda} b_{\mu\lambda} = 0;$$

or, in other words,

$$\sum_{\alpha=1}^{r} c_\alpha \sum_{\lambda=1}^{m} a_{\alpha\lambda} b_{\mu\lambda} + \sum_{x=1}^{s} d_x \sum_{\lambda=1}^{m} x_{x\lambda} b_{\mu\lambda} = 0; \tag{5}$$

The first sum in (5) is, by virtue of (2), equal to $c_\mu$; and the second sum is equal to 0 by virtue of (3); hence, we have $c_\mu = 0$ for $\mu = 1, 2, \ldots, r$.

But then (4) reduces to:

$$\sum_{x=1}^{s} d_x x_{x\lambda} = 0 \quad (\lambda = 1, 2, \ldots, m),$$

and this means that all $d_x = 0$, as the vectors $\mathbf{x}_1, \ldots, \mathbf{x}_s$ are linearly independent.

Now, it is known that the matrix:

$$P = \begin{bmatrix} a_{11} & \cdots & a_{r1} & x_{11} & \cdots & x_{s1} \\ a_{12} & \cdots & a_{r2} & x_{12} & \cdots & x_{s2} \\ \cdots & \cdots & \cdots & \cdots\cdots & \cdots & \cdots \\ a_{1m} & \cdots & a_{rm} & x_{1m} & \cdots & x_{sm} \end{bmatrix}$$

is non-singular when $|P| \neq 0$. We will change the notation, and simplify the notation as follows: $\mathbf{x}_1 = \mathbf{a}_{r+1}, \ldots, \mathbf{x}_s = \mathbf{a}_m$, so:

$$P = (a_{\alpha\beta}) \quad (\alpha, \beta = 1, 2, \ldots, m).$$

We now have the equality:

$$\sum_{\beta=1}^{m} a_{\alpha\beta} y_{\gamma\beta} = \sigma_{\alpha\gamma} \quad (\alpha = 1, 2, \ldots, m; \gamma \text{ is one of } 1, 2, \ldots, m). \tag{6}$$

For a fixed $\gamma$ the equation (6) is a system of $m$ linear equations with $m$ unknowns $y_{\gamma 1}, y_{\gamma 2}, \ldots, y_{\gamma m}$; the determinant of this system is $|P| \neq 0$; therefore, this system (for a fixed $\gamma$) has a unique solution $\mathbf{y}_\gamma = (\mathbf{y}_{\gamma 1}, \mathbf{y}_{\gamma 2}, \ldots, \mathbf{y}_{\gamma m})$. But condition (2) implies that for $\gamma = 1, 2, \ldots, r$ we have $\mathbf{y}_\gamma = \mathbf{b}_\gamma$. We will also, for $\gamma = r + 1, \ldots, m$, denote $\mathbf{y}_\gamma = \mathbf{b}_\gamma$, that is, we determine the values $b_{\alpha\beta}$ also for $\alpha = r + 1, \ldots, m$.

If we now write:

$$Q = \begin{bmatrix} b_{11} & b_{12} & \cdots & b_{1m} \\ b_{21} & b_{22} & \cdots & b_{2m} \\ \cdots & \cdots & \cdots & \cdots \\ b_{m1} & b_{m2} & \cdots & b_{mm} \end{bmatrix}$$

then from (2) we have: $QP = E$, which implies that $Q = P^{-1}$, and, therefore, also a non-singular matrix.

Now, using our notation, we have:

$$K = \begin{bmatrix} a_{11} & \cdots & a_{r1} \\ a_{12} & \cdots & a_{r2} \\ \cdots & \cdots & \cdots \\ a_{1m} & \cdots & a_{rm} \end{bmatrix}, \quad L = \begin{bmatrix} b_{11} & b_{12} & \cdots & b_{1m} \\ \cdots & \cdots & \cdots & \cdots \\ b_{m1} & b_{m2} & \cdots & b_{mm} \end{bmatrix};$$

and hence, by (2), we have:

$$LP = \begin{bmatrix} 1 & 0 & \cdots & 0 & 0 & \cdots & 0 \\ 0 & 1 & \cdots & 0 & 0 & \cdots & 0 \\ \cdots & \cdots & \cdots & \cdots & \cdots & \cdots & \cdots \\ 0 & 0 & \cdots & 1 & 0 & \cdots & 0 \\ 0 & 0 & \cdots & 0 & 0 & \cdots & 0 \\ \cdots & \cdots & \cdots & \cdots & \cdots & \cdots & \cdots \\ 0 & 0 & \cdots & 0 & 0 & \cdots & 0 \end{bmatrix}; P^{-1}K = \left.\begin{bmatrix} 1 & 0 & \cdots & 0 \\ 0 & 1 & \cdots & 0 \\ \cdots & \cdots & \cdots & \cdots \\ 0 & 0 & \cdots & 1 \\ 0 & 0 & \cdots & 0 \\ \cdots & \cdots & \cdots & \cdots \\ 0 & 0 & \cdots & 0 \end{bmatrix}\right\} \begin{matrix} r \\ \\ \\ m-r \end{matrix} \tag{7}$$

Hence we have:

$$P^{-1}IP = (P^{-1}K)(LP) = \begin{bmatrix} 1 & 0 & \cdots & 0 & 0 & \cdots & 0 \\ 0 & 1 & \cdots & 0 & 0 & \cdots & 0 \\ \cdots & \cdots & \cdots & \cdots & \cdots & \cdots & \cdots \\ 0 & 0 & \cdots & 1 & 0 & \cdots & 0 \\ 0 & 0 & \cdots & 0 & 0 & \cdots & 0 \\ \cdots & \cdots & \cdots & \cdots & \cdots & \cdots & \cdots \\ 0 & 0 & \cdots & 0 & 0 & \cdots & 0 \end{bmatrix} = E_m^{(r)}$$

which completes the proof of our theorem.                                                   □

**Example.** Consider the matrix:

$$I = \begin{bmatrix} 3 & 2 & -1 \\ 2 & 3 & -1 \\ 10 & 10 & -4 \end{bmatrix};$$

it is easy to see, that the matrix $I$ is idempotent. We may take:

$$K = \begin{bmatrix} 1 & 0 \\ 1 & -1 \\ 4 & -2 \end{bmatrix}, \quad L = \begin{bmatrix} 3 & 2 & -1 \\ 1 & -1 & 0 \end{bmatrix},$$

and we may now check that

$$KL = I, \quad LK = E_2 = \begin{bmatrix} 1 & 0 \\ 0 & 1 \end{bmatrix}.$$

We now extend the matrix $K$ into a square non-singular matrix $P$ of order 3:

$$P = \begin{bmatrix} 1 & 0 & x \\ 1 & -1 & y \\ 4 & -2 & z \end{bmatrix};$$

and so $x, y, z$ satisfy the equations:

$$3x + 2y - z = 0, \quad x - y = 0.$$

We may set $y = 1$, and hence find: $x = 1, z = 5$. So:

$$P = \begin{bmatrix} 1 & 0 & 1 \\ 1 & -1 & 1 \\ 4 & -2 & 5 \end{bmatrix}.$$

We now extend $L$ into $P^{-1}$:

$$P^{-1} = \begin{bmatrix} 3 & 2 & -1 \\ 1 & -1 & 0 \\ x & y & z \end{bmatrix};$$

where now $x, y, z$ satisfy the conditions:

$$x + yz + 4z = 0, \quad -y - 2z = 0, \quad x + y + 5z = 1;$$

whence:

$$P^{-1} = \begin{bmatrix} 3 & 2 & -1 \\ 1 & -1 & 0 \\ -2 & -2 & 1 \end{bmatrix}.$$

We may now verify that $PP^{-1} = P^{-1}P = E_3$; and $P^{-1}IP = E_3^{(2)} = \begin{bmatrix} 1 & 0 & 0 \\ 0 & 1 & 0 \\ 0 & 0 & 1 \end{bmatrix}.$

**§3.** Let $I, K, L, P$ have the same values as they did in **§2**. Then $P^{-1}IP = E_m^{(r)}$. Next, let $A$ be some matrix of order $m$ and rank $r$, which has $I$ as a two-sided identity. Then (by **§1**) the matrix $LAK = A_1$ is a non-singular matrix of order $r$, and $A = KA_1L$. Let now:

$$A = \begin{bmatrix} a_{11} & a_{12} & \cdots & a_{1r} \\ a_{21} & a_{22} & \cdots & a_{2r} \\ \cdots & \cdots & \cdots & \cdots \\ a_{r1} & a_{r2} & \cdots & a_{rr} \end{bmatrix}$$

(these values $a_{\alpha\beta}$ have nothing to do with the same notation $a_{\alpha\beta}$ used in §2). Then (by (7) in §2):

$$P^{-1}AP = (P^{-1}K)A_1(LP) = \begin{bmatrix} a_{11} & \cdots & a_{1r} & 0 & \cdots & 0 \\ \cdots & \cdots & \cdots & \cdots & \cdots & \cdots \\ a_{r1} & \cdots & a_{rr} & 0 & \cdots & 0 \\ 0 & \cdots & 0 & 0 & \cdots & 0 \\ \cdots & \cdots & \cdots & \cdots & \cdots & \cdots \\ 0 & \cdots & 0 & 0 & \cdots & 0 \end{bmatrix}.$$

Hence, the matrix $P^{-1}AP$ has the following form: it is the matrix $A_1$, but surrounded from below and on the right side by $m - r$ rows of zeroes. We will abbreviate this as:

$$A_m^{(r)} = P^{-1}AP = \begin{bmatrix} A_1 & 0_{r,s} \\ 0_{s,r} & 0_s \end{bmatrix} \tag{8}$$

where $s = m - r$.

If $\mathfrak{G}_1$ denotes the ordinary group of all matrices of order $m$ and which has a two-sided identity $I$, then $P^{-1}\mathfrak{G}_1 P$ is the group, simply isomorphic to $\mathfrak{G}_1$, of all matrices of order $m$ and of type (8), where $A_1$ ranges over all possible non-singular matrices of order $r$.

**§4.** We will now consider the problem: find all the right identities $X = (x_{\alpha\beta})$ for the matrix $A_m^{(r)}$ (see (8), §3). The condition $A_m^{(r)} X = A_m^{(r)}$ gives:

$$\sum_{\lambda=1}^{r} a_{\alpha\lambda} x_{\lambda\beta} = a_{\alpha\beta} \qquad (\alpha, \beta = 1, 2, \ldots, r). \tag{9}$$

For a fixed $\beta$, (9) is a system of linear equations with $r$ unknowns $x_{1\beta}, x_{2\beta}, \ldots, x_{r\beta}$; the determinant of this system is $|A_1| \neq 0$; so this system always has exactly one solution; therefore, this unique solution is: $x_{\beta\beta} = 1$, $x_{\lambda\beta} = 0$ for $\lambda \neq \beta, \lambda \leq r$. If $\beta = r + 1, \ldots, m$, then the right-hand sides of (9) are all $= 0$, and there all $x_{\lambda\beta} = 0$ (for $\beta > r$). Finally, for $\lambda = r + 1, \ldots, m$ all $x_{\lambda\beta}$ remain arbitrary. Hence:

$$X = \begin{bmatrix} 1 & 0 & \cdots & 0 & 0 & \cdots & 0 \\ 0 & 1 & \cdots & 0 & 0 & \cdots & 0 \\ \cdots & \cdots & \cdots & \cdots & \cdots & \cdots & \cdots \\ 0 & 0 & \cdots & 1 & 0 & \cdots & 0 \\ \alpha_{11} & \alpha_{12} & \cdots & \alpha_{1r} & \alpha_{1,r+1} & \cdots & \alpha_{1m} \\ \cdots & \cdots & \cdots & \cdots & \cdots & \cdots & \cdots \\ \alpha_{s1} & \alpha_{s2} & \cdots & \alpha_{sr} & \alpha_{s,r+1} & \cdots & \alpha_{sm} \end{bmatrix} \tag{10}$$

where $s = m - r$, and all $\alpha_{r\lambda}$ are arbitrary. It is easy to check that every matrix of type (10) is a right identity for $A_m^{(r)}$ and, more generally, for every matrix of type (8), in particular for $E_m^{(r)}$.

Analogously, we find that every left identity for matrices of the type (8) (and in particular for $E_m^{(r)}$) has the form:

$$Y = \begin{bmatrix} 1 & 0 & \cdots & 0 & \beta_{11} & \cdots & \beta_{1s} \\ 0 & 1 & \cdots & 0 & \beta_{21} & \cdots & \beta_{2s} \\ \cdots & \cdots & \cdots & \cdots & \cdots & & \cdots \\ 0 & 0 & \cdots & 1 & \beta_{r1} & \cdots & \beta_{rs} \\ 0 & 0 & \cdots & 0 & \beta_{r+1,1} & \cdots & \beta_{r+1,s} \\ \cdots & \cdots & \cdots & \cdots & \cdots & & \cdots \\ 0 & 0 & \cdots & 0 & \beta_{m1} & \cdots & \beta_{ms} \end{bmatrix} \tag{10a}$$

where all $\beta_{x\lambda}$ are arbitrary. We see that matrices of the type $Y$ are transposes ($X'$) of matrices of type $X$.

Let then $Z$ be a two-sided identity for $A_m^{(r)}$; then $Z$ is both of the form (10) as well as (10a); hence, $Z$ will have the form:

$$Z = \begin{bmatrix} 1 & 0 & \cdots & 0 & 0 & \cdots & 0 \\ 0 & 1 & \cdots & 0 & 0 & \cdots & 0 \\ \cdots & \cdots & \cdots & \cdots & \cdots & \cdots & \cdots \\ 0 & 0 & \cdots & 1 & 0 & \cdots & 0 \\ 0 & 0 & \cdots & 0 & \gamma_{11} & \cdots & \gamma_{1s} \\ \cdots & \cdots & \cdots & \cdots & \cdots & \cdots & \cdots \\ 0 & 0 & \cdots & 0 & \gamma_{s1} & \cdots & \gamma_{ss} \end{bmatrix} \tag{10b}$$

It is easy to see that the ranks of the matrices $X, Y, Z$ are never less than $r$; but they can be bigger; these matrices can be both non-singular and non-regular; if their rank $> r$, then they can be non-idempotent. But the group property for matrices of each of these three types is fulfilled, i.e. the product of two right (or left, or two-sided) identities is again a right (or left, or two-sided) identity. If we also set the condition that the right (or left, or two-sided) identity for $A_m^{(r)}$ have rank $r$, then we can easily show (by using Kronecker's famous Theorem) that: any right identity of rank $r$ must have the form:

$$X_r = \begin{bmatrix} 1 & 0 & \cdots & 0 & 0 & \cdots & 0 \\ 0 & 1 & \cdots & 0 & 0 & \cdots & 0 \\ \cdots & \cdots & \cdots & \cdots & \cdots & \cdots & \cdots \\ 0 & 0 & \cdots & 1 & 0 & \cdots & 0 \\ \alpha_{11} & \alpha_{12} & \cdots & \alpha_{1r} & 0 & \cdots & 0 \\ \cdots & \cdots & \cdots & \cdots & \cdots & \cdots & \cdots \\ \alpha_{s1} & \alpha_{s2} & \cdots & \alpha_{sr} & 0 & \cdots & 0 \end{bmatrix} \tag{10c}$$

and any left identity of rank $r$ must have the form:

$$Y = \begin{bmatrix} 1 & 0 & \cdots & 0 & \beta_{11} & \cdots & \beta_{1s} \\ 0 & 1 & \cdots & 0 & \beta_{21} & \cdots & \beta_{2s} \\ \cdots & \cdots & \cdots & \cdots & \cdots & \cdots & \cdots \\ 0 & 0 & \cdots & 1 & \beta_{r1} & \cdots & \beta_{rs} \\ 0 & 0 & \cdots & 0 & 0 & \cdots & 0 \\ \cdots & \cdots & \cdots & \cdots & \cdots & \cdots & \cdots \\ 0 & 0 & \cdots & 0 & 0 & \cdots & 0 \end{bmatrix} \tag{10d}$$

(where $\alpha_{x\lambda}, \beta_{x\lambda}$ are arbitrary). Hence, there is only a single two-sided identity for $A_m^{(r)}$ of rank $r \in E_m^{(r)}$. The matrices $X_r$ and $Y_r$ are idempotent (and hence regular).

The preceding pages allows for finding all right, left, and two-sided identity elements for an arbitrary regular matrix of order $m$ and rank $r < m$. Hence, we have:

**Theorem.** *All ordinary matrices of order $m$ and rank $r < m$, which share a two-sided identity $I$ of order $m$ and ran $r$, all share the same set of right identities, the same set of left identities, and the same set of two-sided identities. More precisely, taking $P^{-1}IP = E_m^{(r)}$, then the right identities will all have the form: $PXP^{-1}$, of which those of rank $r$ will have the form: $PX_rP^{-1}$; the left identities will have the form: $PYP^{-1}$, of which those of rank $r$ will have the form: $PY_rP^{-1}$; and finally, the two-sided identities will have the form: $PZP^{-1}$; where $X, Y, X_r, Y_r, Z$ are given by formulas (10)–(10d).*

**Example.** Let us find all right, left, and two-sided identities for those matrices of order 3 and rank 2 which all have as a two-sided identity the matrix $I = \begin{bmatrix} 3 & 2 & -1 \\ 2 & 3 & -1 \\ 10 & 10 & -4 \end{bmatrix}$ (see the example in §2). We have:

$$P = \begin{bmatrix} 1 & 0 & 1 \\ 1 & -1 & 1 \\ 4 & -2 & 5 \end{bmatrix}, \qquad P^{-1} = \begin{bmatrix} 3 & 2 & -1 \\ 1 & -1 & 0 \\ -2 & -2 & 1 \end{bmatrix}$$

So, we will find all right identities have the form:

$$PXP^{-1} = P \begin{bmatrix} 1 & 0 & 0 \\ 0 & 1 & 0 \\ \alpha_1 & \alpha_2 & \alpha_3 \end{bmatrix} P^{-1} =$$

$$= \begin{bmatrix} 3+3\alpha_1+\alpha_2-2\alpha_3 & 2+2\alpha_1-\alpha_2-2\alpha_3 & -1-\alpha_1+\alpha_3 \\ 2+3\alpha_1+\alpha_2-2\alpha_3 & 3+2\alpha_1-\alpha_2-2\alpha_3 & -1-\alpha_1+\alpha_3 \\ 10+15\alpha_1+5\alpha_2-10\alpha_3 & 10+10\alpha_1-5\alpha_2-10\alpha_3 & -4-5\alpha_1+5\alpha_3 \end{bmatrix};$$

of these, the ones with rank 2 are those with $\alpha_3 = 0$:

$$PX_2P^{-1} = P \begin{bmatrix} 1 & 0 & 0 \\ 0 & 1 & 0 \\ \alpha_1 & \alpha_2 & 0 \end{bmatrix} P^{-1} = \begin{bmatrix} 3+3\alpha_1+\alpha_2 & 2+2\alpha_1-\alpha_2 & -1-\alpha_1 \\ 2+3\alpha_1+\alpha_2 & 3+2\alpha_1-\alpha_2 & -1-\alpha_1 \\ 10+15\alpha_1+5\alpha_2 & 10+10\alpha_1-5\alpha_2 & -4-5\alpha_1 \end{bmatrix};$$

All the left identities will have the form:

$$PYP^{-1} = P \begin{bmatrix} 1 & 0 & \beta_1 \\ 0 & 1 & \beta_2 \\ 0 & 0 & \beta_3 \end{bmatrix} =$$

$$= \begin{bmatrix} 3-2\beta_1-2\beta_3 & 2-2\beta_1-2\beta_3 & -1+\beta_1+\beta_3 \\ 2-2\beta_1+2\beta_2-2\beta_3 & 3-2\beta_1+2\beta_2-2\beta_3 & -1+\beta_1-\beta_2+\beta_3 \\ 10-8\beta_1+4\beta_2-10\beta_3 & 10-8\beta_1+4\beta_2-10\beta_3 & -4+4\beta_1-2\beta_2+5\beta_3 \end{bmatrix};$$

of these, the ones with rank 2 are those with $\beta_3 = 0$:

$$PY_2P^{-1} = P \begin{bmatrix} 1 & 0 & \beta_1 \\ 0 & 1 & \beta_2 \\ 0 & 0 & 0 \end{bmatrix} = \begin{bmatrix} 3-2\beta_1-2 & 2-2\beta_1 & -1+\beta_1 \\ 2-2\beta_1+2\beta_2 & 3-2\beta_1+2\beta_2 & -1+\beta_1-\beta_2 \\ 10-8\beta_1+4\beta_2 & 10-8\beta_1+4\beta_2 & -4+4\beta_1-2\beta_2 \end{bmatrix};$$

Finally, we obtain the two-sided identities by setting, for the right identities,

$$\alpha_1 = \alpha_2 = 0, \alpha_3 = \gamma, \quad \text{or setting, for the left, } \beta_1 = \beta_2 = 0, \beta_3 = \gamma:$$

$$PZP^{-1} = P \begin{bmatrix} 1 & 0 & 0 \\ 0 & 1 & 0 \\ 0 & 0 & \gamma \end{bmatrix} P^{-1} = \begin{bmatrix} 3 - 2\gamma & \gamma 2 - 2\gamma & -1 + \gamma \\ 2 - 2\gamma & 3 - 2\gamma & -1 + \gamma \\ 10 - 10\gamma & 10 - 10\gamma & -4 + 5\gamma \end{bmatrix}.$$

For $\gamma = 0$, we get $I$, as expected; and for $\gamma = 1$ we get $E_3$.

We note further that this two-sided identity $PZP^{-1}$ is represented in the form

$PZP^{-1} = I + \gamma N$, where $N = \begin{bmatrix} -2 & -2 & 1 \\ -2 & -2 & 1 \\ -10 & -10 & 5 \end{bmatrix}$ is a matrix such that

$$IN = 0.$$

**§5.** Let us now post the "inverse" problem: find all matrices (of order $m$) $A = (a_{\alpha\beta})$, for which $E_m^{(r)}$ is a right identity.

Let us denote temporarily: $E_m^{(r)} = (x_{\alpha\beta})$; then the condition $AE_m^{(r)} = A$ becomes:

$$\sum_{\lambda=1}^{m} a_{\alpha\lambda} x_{\lambda\beta} = a_{\alpha\beta} \quad (\alpha, \beta = 1, 2, \ldots, m). \tag{11}$$

But we know that $x_{\alpha\beta} = 0$ for $\alpha \neq \beta$; thus, (11) reduces to

$$a_{\alpha\beta} x_{\beta\beta} = a_{\alpha\beta} \quad (\alpha, \beta = 1, 2, \ldots, m). \tag{11a}$$

In addition, we have $x_{\beta\beta} = 1$ for $\beta = 1, 2, \ldots, r$, and for $\beta = r + 1, \ldots, m$ we have $x_{\beta\beta} = 0$. Hence, for $\beta \leq r$ (11a) is identically true, i.e. $a_{\alpha\beta}$ can be chosen arbitrarily for $\beta \leq r$; and for $\beta > r$ equation (11a) yields: $a_{\alpha\beta} = 0$.

We hence have that the matrix $A$ can be written in the form:

$$A = \begin{bmatrix} \alpha_{11} & \cdots & \alpha_{1r} & 0 & \cdots & 0 \\ \alpha_{21} & \cdots & \alpha_{2r} & 0 & \cdots & 0 \\ \cdots & \cdots & \cdots & \cdots & \cdots & \cdots \\ \alpha_{m1} & \cdots & \alpha_{mr} & 0 & \cdots & 0 \end{bmatrix} \tag{12}$$

where the $a_{\alpha\beta}$ (for $\alpha = 1, \ldots, m; \beta = 1, \ldots, r$) are arbitrary.

Analogously, we may show, that any matrix $B$ (of order $m$) for which $E_m^{(r)}$ is a left identity, has the form:

$$B = \begin{bmatrix} b_{11} & b_{12} & \cdots & b_{1m} \\ \cdots & \cdots & \cdots & \cdots \\ b_{r1} & b_{r2} & \cdots & b_{rm} \\ 0 & 0 & \cdots & 0 \\ \cdots & \cdots & \cdots & \cdots \\ 0 & 0 & \cdots & 0 \end{bmatrix} \tag{12b}$$

where the $b_{\alpha\beta}$ (for $\alpha = 1, \ldots, r; \beta = 1, \ldots, m$) are arbitrary.

It is easy to see that the group property is fulfilled for the set of matrices of the form (12), i.e. this set is a generalized group; the same applies to the set of matrices (12a). The ranks of these matrices (12) and (12a) is no greater than $r$. Note that the matrix $X_r$ (see §4, (10c)) is of type (12), and the matrix $Y_r$ (see §4, (10d)) is of type (12a), that is:

$$X_r E_m^{(r)} = X_r, \quad E_m^{(r)} Y_r = Y_r.$$

The intersection of the sets of matrices (12) and (12a) is a matrix of the form:

$$C = \begin{bmatrix} c_{11} & \cdots & c_{1r} & 0 & \cdots & 0 \\ \cdots & \cdots & \cdots & \cdots & \cdots & \cdots \\ c_{r1} & \cdots & c_{rr} & 0 & \cdots & 0 \\ 0 & \cdots & 0 & 0 & \cdots & 0 \\ \cdots & \cdots & \cdots & \cdots & \cdots & \cdots \\ 0 & \cdots & 0 & 0 & \cdots & 0 \end{bmatrix} \qquad (12\text{B})$$

For these, and only these, matrices $C$ the matrix $E_m^{(r)}$ is a two-sided identity. This group includes those matrices $A_m^{(r)}$ considered in §3 (see (8)); but those matrices were of rank $r$, whereas the rank of our matrices $C$ can be $\leq r$. Note that such matrices $C$, if they are of rank $r$, are necessarily regular.

# GENERALISED GROUPS OF SOME
# TYPES OF INFINITE MATRICES

## A. K. Sushkevich (Kharkiv)

**§1.** In this article, I will consider infinite bounded matrices of the form:

$$A = \begin{pmatrix} a_{11} & a_{12} & \cdots \\ a_{21} & a_{23} & \cdots \\ \cdots & \cdots & \cdots \end{pmatrix} \tag{1}$$

where the coefficients $a_{\alpha\beta}$ are all real numbers, and the vectors $\mathbf{a}_\lambda = (a_{\lambda 1}, a_{\lambda 2}, \dots)$ for $\lambda = 1, 2, \dots$ *ad inf.*, form an orthonormal system of vectors in Hilbert space. This means:

$$\sum_{\lambda=1}^{\infty} a_{\alpha\lambda} a_{\beta\lambda} = \delta_{\alpha\beta} = \begin{cases} 1 & \text{if } \alpha = \beta, \\ 0 & \text{if } \alpha \neq \beta. \end{cases} \tag{2}$$

If the system of vectors $\mathbf{a}_\lambda$ is a *complete* orthonormal system, then the matrix $A$ is, as usual, called *orthogonal*, and then also

$$\sum_{\lambda=1}^{\infty} a_{\lambda\alpha} a_{\lambda\beta} = \delta_{\alpha\beta} = \begin{cases} 1 & \text{if } \alpha = \beta, \\ 0 & \text{if } \alpha \neq \beta. \end{cases} \tag{3}$$

In other words, the transposed matrix $A'$ is also orthogonal, and $AA' = A'A = E$, the identity matrix, and $A' = A^{-1}$ is an *inverse* matrix for $A$ (which, therefore, exists).

Bearing in mind that the matrix $E$ is also orthogonal, that the product of orthogonal matrices is again orthogonal, and that the associative law holds for the multiplication of bounded matrices, we can draw the following conclusion:

> *The collection of all infinite orthogonal matrices forms an ordinary group.*

**§2.** We will, however, consider such matrices $A$ for which (2) is fulfilled but (3) is not fulfilled, that is, the orthonormal system of vectors $\mathbf{a}_\lambda$ is not complete. Let us call such a matrix $A$ *left semi-orthogonal*. We will say that a left semi-orthogonal matrix $A$ is *proper*, if an infinite number of vectors are required to complete the system of vectors $\mathbf{a}_\lambda$ into a complete one.

**Theorem.** *The collection of all proper left semi-orthogonal matrices forms a generalized associative group, which satisfies the left-sided law of unique reversibility and the right-sided law of unlimited reversibility.*

*Proof.* First, we prove that group property holds for the collection of matrices of this type. Let $A$ and $B = (b_{\alpha\beta})$ be two such matrices, and denote their product $AB = C = (c_{\alpha\beta})$. We will prove that $C$ is a proper left semi-orthogonal matrix. We have

$$c_{\alpha\beta} = \sum_{\lambda=1}^{\infty} a_{\alpha\lambda} b_{\lambda\beta}.$$

All the series that we consider are absolutely convergent (because our matrices are bounded): We have

$$\sum_{\mu=1}^{\infty} c_{\alpha\mu} c_{\beta\mu} = \sum_{\mu=1}^{\infty} \left(\sum_{x=1}^{\infty} a_{\alpha x} b_{x\mu}\right) \left(\sum_{\lambda=1}^{\infty} a_{\alpha\lambda} b_{\lambda\mu}\right) = \sum_{x,\lambda=1} \left(a_{\alpha x} a_{\beta\lambda} \sum_{\mu=1}^{\infty} b_{x\mu} b_{\lambda\mu}\right) =$$

$$= \sum_{x,\lambda=1}^{\infty} a_{\alpha x} a_{\alpha\lambda} \delta_{x\lambda} = \sum_{\lambda=1}^{\infty} a_{x\lambda} a_{\beta\lambda} = \delta_{\alpha\beta}.$$

Hence equation (3) has been shown to hold for the matrix $C$.

Let now $\mathbf{a}' = (a_1', a_2', \dots)$ be any vector, orthogonal to all $a_\lambda$. We set

$$c_\lambda' = \sum_{x=1}^{\infty} a_x' b_{x\lambda}$$

and we will now show that the vector $\mathbf{c}' = (c_1', c_2', c_3', \dots)$ is different from all vectors $\mathbf{c}_\alpha = (c_{\alpha_1}, c_{\alpha_2}, \dots)$, and that it is orthogonal to every $\mathbf{c}_\alpha$.

Suppose now that $\mathbf{c}' = \mathbf{c}_\alpha$, i.e. $c_\lambda' = c_{\alpha\lambda}$ for every $\lambda$; then we have

$$\sum_{x=1}^{\infty} (a_x' - a_{\alpha x}) b_{x\lambda} = 0 \quad (\lambda = 1, 2, \dots, \text{ ad inf.}). \tag{4}$$

Let us multiply (4) by $b_{\mu\lambda}$, and sum with respect to $\lambda$:

$$0 = \sum_{x=1}^{\infty} (a_x' - a_{\alpha x}) \sum_{\lambda=1}^{\infty} b_{\alpha\lambda} b_{\mu\lambda} = \sum_{x=1}^{\infty} (a_x' - a_{\alpha x}) \delta_{x\mu} = a_\mu' - a_{\alpha\mu}.$$

Therefore $a_\mu' = a_{\alpha\mu}$ for every $\mu$, whence $\mathbf{a}' = \mathbf{a}_\alpha$. Thus, as $\mathbf{a}'$ is different from $\mathbf{a}_\alpha$, this is a contradiction, so we must have that $\mathbf{c}'$ is different from $\mathbf{c}_\alpha$. Then we have:

$$\sum_{\mu=1}^{\infty} c_{\alpha\mu} c_\mu' = \sum_{\mu=1}^{\infty} \left(\sum_{x=1}^{\infty} a_{\alpha x} b_{x\mu}\right) \left(\sum_{\lambda=1}^{\infty} a_\lambda' b_{\lambda\mu}\right) =$$

$$= \sum_{x,\lambda=1}^{\infty} \left(a_{\alpha x} a_\lambda' \sum_{\mu=1}^{\infty} b_{\alpha\mu} b_{\lambda\mu}\right) = \sum_{x,\lambda=1}^{\infty} a_{\alpha x} a_\lambda' \delta_{x\lambda} = \sum_{\lambda=1}^{\infty} a_{\alpha\lambda} a_\lambda' = 0.$$

Therefore, the vector $\mathbf{c}'$ is orthogonal to every $\mathbf{c}_\alpha$. If the vector $\mathbf{a}'$ is normalized, then $\mathbf{c}'$ is also normalized, as:

$$\sum_{\alpha=1}^{\infty} c_\alpha'^2 = \sum_{\alpha=1}^{\infty} \left(\sum_{x=1}^{\infty} a_x' b_{x\alpha}\right) \left(\sum_{\lambda=1}^{\infty} a_\lambda' b_{\lambda\alpha}\right) = \sum_{x,\lambda=1}^{\infty} \left(a_x' a_\lambda' \sum_{\alpha=1}^{\infty} b_{x\alpha} b_{\lambda\alpha}\right) = \sum_{\lambda=1}^{\infty} a_\lambda'^2 = 1.$$

We now add the vector $\mathbf{a}'$ to the system of vectors $\mathbf{a}_\alpha$, and at the same time add the vector $\mathbf{c}'$ to the system of vectors $\mathbf{c}_\alpha$. One can now find a normalized vector $\mathbf{a}''$, orthogonal to all $\mathbf{a}_\alpha$ and to $\mathbf{a}'$, and, correspondingly, a normalized vector $\mathbf{c}''$, orthogonal to all $\mathbf{c}_\alpha$ and to $\mathbf{c}'$, etc. without end. Hence, since the matrix $A$ is a proper left semi-orthogonal matrix, it follows that $C$ is also a proper left semi-orthogonal matrix.

We now prove that the left-sided law of unique reversibility holds for the matrices under consideration. Suppose that

$$BA = CA. \tag{5}$$

From equation (2) it follows that $AA' = E$ (the identity matrix). Hence, multiplying both sides of (5) on the right by $A'$, we get $BE = CE$, i.e. $B = C$.

We will now prove that the right-sided law of unlimited reversibility holds. Suppose we have the equality

$$AX = B, \tag{6}$$

for some $A = (a_{\alpha\beta})$, $B = (b_{\alpha\beta})$ of the given type, and $X = (x_{\alpha\beta})$ is an unknown matrix. The equality in (6) yields:

$$\sum_{x=1}^{\infty} a_{\alpha x} x_{x\beta} = b_{\alpha\beta} \quad (\alpha = 1, 2, \ldots, \text{ad inf. and } \beta = 1, 2, \ldots, \text{ad inf.}). \tag{6a}$$

For a fixed $\beta$ this is a system of an infinite number of linear equalities with infinitely many unknowns $x_{1\beta}, x_{2\beta}, \ldots$. This system is orthonormal. Its general solution can be found as follows.

Let $\mathbf{a}'_\rho = (a'_{\rho 1}, a'_{\rho 2}, \ldots)$, where $\rho = 1, 2, \ldots$ ad inf., be a vector which complements the system $\mathbf{a}_\alpha$ into a complete system. Let us, in addition to (6a), take the equation

$$\sum_{\mu=1}^{\infty} a'_{\rho\mu} x_{\mu\beta} = v_{\rho\beta}, \tag{7}$$

where $(v_{1\beta}, v_{2\beta}, \ldots)$ is an *arbitrary* vector. Then the general solution to the equations in (6a) will be given as:

$$x_{\alpha\beta} = \sum_{\nu=1}^{\infty} a_{\nu x} b_{\nu\beta} + \sum_{\sigma=1}^{\infty} a'_{\sigma x} v_{\sigma\beta}. \tag{8}$$

This formula defines a matrix $X$, and hence can do so in infinitely many ways. But we still need to prove that any such matrix $X$ can be chosen as a proper left semi-orthogonal matrix. For this, we choose $v_{\alpha\beta}$ such that

$$\mathbf{v}_\sigma = (v_{\sigma 1}, v_{\sigma 2}, \ldots), \quad (\sigma = 1, 2, \ldots, \text{ ad inf.})$$

will be normalized vectors, mutually orthogonal and orthogonal to all vectors $\mathbf{b}_\alpha = (b_{\alpha 1}, b_{\alpha 2}, \ldots)$.

We can, for example, choose as $\mathbf{v}_\sigma$ *all* vectors, which complement the system $\mathbf{b}_\alpha$ into a new one, or only *some* (but infinitely many, of course) of these vectors. Let us prove that then condition (2) is fulfilled. We have:

$$\sum_{\beta=1}^{\infty} x_{\alpha\beta} x_{\lambda\beta} = \sum_{\mu,\nu=1}^{\infty} \sum_{\beta=1}^{\infty} a_{\mu x} a_{\nu\lambda} b_{\nu\beta} b_{\nu\beta} + \sum_{\nu,\sigma=1}^{\infty} \sum_{\beta=1}^{\infty} a_{\mu x} a'_{\sigma\lambda} b_{\mu\beta} v_{\sigma\beta} +$$
$$+ \sum_{\rho,\nu=1}^{\infty} \sum_{\beta=1}^{\infty} a_{\nu\lambda} a'_{\rho x} b_{\nu\beta} v_{\rho\beta} + \sum_{\rho,\sigma=1}^{\infty} \sum_{\beta=1}^{\infty} a'_{\rho x} a'_{\sigma\lambda} v_{\rho\beta} v_{\sigma\beta} =$$
$$= \sum_{\mu=1}^{\infty} a_{\mu x} a_{\mu\lambda} + \sum_{\rho=1}^{\infty} a'_{\rho x} a'_{\rho\lambda} = \delta_{x\lambda},$$

because this is equation (3) for the complete system of orthogonal vectors $\mathbf{a}_\mu$ and $\mathbf{a}'_\beta$. Additionally, we have

$$\sum_{\beta=1}^{\infty} b_{\mu\beta} b_{\nu\beta} = \delta_{\mu\nu}, \quad \sum_{\beta=1}^{\infty} b_{\mu\beta} v_{\sigma\beta} = \sum_{\beta=1}^{\infty} b_{\nu\beta} v_{\rho\beta} = 0, \quad \sum_{\beta=1}^{\infty} v_{\rho\beta} v_{\sigma\beta} = \delta_{\rho\sigma}.$$

This proves that the vectors $\mathbf{x}_\lambda = (x_{\lambda 1}, x_{\lambda 2}, \dots)$ forms an orthonormal system. We will now uncover when this system will be complete. For that, equation (3) must be satisfied, that is:

$$\delta_{\beta\gamma} = \sum_{x=1}^\infty x_{x\beta} x_{x\gamma} = \sum_{\mu,\nu=1}^\infty \sum_{x=1}^\infty a_{\mu x} a_{\nu x} b_{\mu\beta} b_{\nu\gamma} + \sum_{\mu,\sigma=1}^\infty \sum_{x=1}^\infty a_{\mu x} a'_{\sigma x} b_{\mu\beta} v_{\sigma\gamma} +$$

$$+ \sum_{\rho,\sigma=1}^\infty \sum_{x=1}^\infty a'_{\rho x} a'_{\sigma x} v_{\rho\beta} v_{\sigma\beta} = \sum_{\mu=1}^\infty b_{\mu\beta} b_{\mu\gamma} + \sum_{\rho=1}^\infty v_{\rho\beta} v_{\rho\gamma}$$

(which holds because the middle two sums are equal to zero). Hence, the orthonormal system of vectors $(\mathbf{b}_1, \mathbf{b}_2, \dots, \mathbf{v}_1, \mathbf{v}_2, \dots)$ is complete.

Conversely, if this system is complete, then:

$$\sum_{x=1}^\infty x_{x\beta} x_{x\gamma} = \delta_{\beta\gamma},$$

that is, the vectors $\mathbf{x}_\lambda$ form a complete orthonormal system, and the matrix $X$ is orthogonal.

However, we can, in an infinite number of ways choose, from the latter half of the system $(\mathbf{b}_1, \mathbf{b}_2, \dots, \mathbf{v}_1, \mathbf{v}_2, \dots)$, many vectors $v_{\alpha_1}, v_{\alpha_2} \dots$ such that the rest of the vectors $\mathbf{v}_\lambda$ remains infinite. We will prove that in this case, the vectors $\mathbf{x}_\lambda$, defined by equation (8), will form an incomplete orthonormal system, with the property that an infinite number of vectors will be required to complement it to a complete system. For this, we need to prove that there is an infinite number of vectors orthogonal to all $\mathbf{x}_\lambda$. But we know, that (according to our condition) there are infinitely many vectors orthogonal to each other, to all $\mathbf{b}_\lambda$ and to all *chosen* $v_\lambda$ (we will henceforth denote these chosen vectors $v_\lambda$ as $v_1, v_2, \dots$ and not as $v_{\alpha_1}, v_{\alpha_2}, \dots$).

So, let $\mathbf{y} = (y_1, y_2, \dots)$ be one of the infinitely many vectors orthogonal to all $\mathbf{b}_\lambda$ and the chosen $v_\lambda$. This proves that the matrix $X$ in equation (6) can be defined in an infinite number of ways as a proper left semi-orthogonal matrix. This completes the proof of the theorem. □

**§3.** All the conclusions about proper left semi-orthogonal matrices given in §2 can be "inverted" if in all our matrices we rearrange the rows into columns. The matrices obtained in this way are called *right semi-orthogonal*, and similarly *proper*. It is now easy to see that the following theorem holds:

**Theorem.** *The collection of all proper right semi-orthogonal matrices forms a generalized associative group, which satisfies the right-sided law of unique reversibility and the left-sided law of unlimited reversibility.*

**§4.** As an example, let us consider the following special case of proper left semi-orthogonal matrices which form a generalized group of the type considered in §2:

Let the matrix $A$ be composed as follows:

(1) each row in $A$ consists of zeroes, except for one entry in the row, which is equal to 1;

(2) each column in $A$ consists of zeroes, except for *at most* one entry in the column, which is equal to 1;

(3) there are infinitely many columns in $A$ consisting only of zeroes.

We will prove that the collection of all matrices of this type forms a generalized group of the same type as in §2. The associative law and the left-sided law of unique reversibility were proved already in §2. We will first prove that the group property holds for these matrices. Let $A$ and $B$ be two such matrices, with $AB = C$. Consider the $\alpha$-th row of the matrix $C$; it consists of the composition of the $\alpha$-th row of the matrix $A$ with all the columns of the matrix $B$. Suppose the element 1 in the $\alpha$-th row of the matrix $A$ lies in the $\beta$-th column. If there is no unit at all in a given column of $B$ (that is, that column consists of zeroes), or if the unit is not in the $\beta$-th place, then the corresponding element in the $\alpha$-th row of the matrix $C$ will be 0. But due to the fact that the $\beta$-th row of $B$ necessarily has one and only one element 1, – for example, in the $\gamma$-th column – then the composition of the $\alpha$-th row of the matrix $A$ with the $\gamma$-th column of the matrix $B$ will be identical as a result of the composition, that is, in the $\alpha$-th row of the matrix $C$, the $\gamma$-th column will have the element 1, and all other elements of the $\alpha$-th row will be 0. Hence property 1) holds for the matrix $C$.

Let us now consider the $\beta$-th column of the matrix $C$. It is form by the composition of all the rows of the matrix $A$ with the $\beta$-th column of the matrix $B$. But, if the $\beta$-th column of $B$ is "empty" (that is, it consists of only zeroes), then also in $C$ the $\beta$-th row will be empty. This proves properties 2) and 3) for $C$. At the same time, there will be "more" empty columns in $C$ than in $B$.

We now turn to the proof of the right-sided law of unrestricted reversibility for matrices of our type. Let $A$ and $B$ be two matrices satisfying properties 1), 2), and 3). We will prove that we can find infinitely many matrices $X$, of the same type, such that $AX = B$.

Suppose that the $\alpha$-th row of $A$ has the element 1 in the $\beta$-th column, and that the $\alpha$-th row of $B$ has an element 1 in the $\gamma$-th column. Then in $X$, in the $\beta$-th row we will put the element 1 in the $\gamma$-th column, and all other elements in the $\beta$-th row of $X$ with be taken as 0. But this does not determine all rows in $X$: infinitely many of them are still undefined. We will fill these undefined rows in $X$ with zeroes, putting a unit in only one place in each row, but in such a way that the properties 2) and 3) remain true for $X$ (this is possible, because the numbers $\gamma$ do not exhaust the set of natural numbers; there are infinitely many natural numbers not equal to any $\gamma$).

This proves our statement about matrices with the properties 1), 2), and 3).

By swapping rows for columns in all our matrices, we obtain a another type of proper right semi-orthogonal matrices, which constitute a generalized group of type §3.

We remark that the matrices with properties 1), 2), and 3) of this section are nothing but the so-called "insufficient" permutations on a countable set of symbols (objects), i.e. permutations in which all the symbols appear, and without duplicates, in the lower row. At the same time, the proper matrices correspond to such permutations in which there is no infinite number of symbols in their lower row.

We also remark that the "excessive" permutations of a countable set of symbols (i.e. those in which all the symbols appear in the bottom row, and the same one can occur several times) do not give us semi-orthogonal matrices.

I will consider the different types of permutations on a countable set of symbols elsewhere.



%\vspace{0.5cm}

\end{document}